\documentclass[preprint]{elsarticle}
\usepackage[english]{babel}
\usepackage{graphics}
\usepackage{graphicx}
\usepackage{epsfig}
\usepackage{amssymb}
\usepackage{amsthm}
\usepackage{amsmath}
\usepackage{amsfonts}
\usepackage{pstricks,pst-plot}
\usepackage{hyperref}
\usepackage{epstopdf}
\usepackage{subfigure}
\usepackage{amssymb}
\usepackage{dsfont}

\newcommand{\pippo}{\mbox{\emph {sum}}}

\newcommand{\roots}{\mbox{\emph {roots}}}
\newcommand{\trace}{\mbox{\rm {trace}}}
\newcommand{\ones}{\mbox{\emph {ones}}}
\newcommand{\zeros}{\mbox{\emph {zeros}}}
\newcommand{\planerot}{\mbox{\emph {planerot}}}
\newcommand{\diag}{\mbox{\rm {diag}}}
\newcommand{\spec}{\mbox{\rm {spec}}}
\newcommand{\B}[1]{\mbox{\boldmath $#1$}}
\newtheorem{lemma}{Lemma}

\newtheorem{remark}{Remark}
\newtheorem{proposition}{Proposition}

\newenvironment{code}{%
                           \mathcode`\:="603A  
                           \def\colon{\mathchar"303A}
                           \par
                           \upshape
                           \begin{list} 
                              {} {\leftmargin = 0.0cm}
                           \item[]
                           \begin{tabbing}
                              \hspace*{.3in} \= \hspace*{.3in} \=
                              \hspace*{.3in} \= \hspace*{.3in} \=

                          \hspace*{.3in} \= \hspace*{.3in} \= \kill
                          }{\end{tabbing}\end{list}}
\input epsf.tex
\def\fig#1#2#3{
\vskip #3pc \centerline{{\bf Fig.~#1.} {\it #2}}  }

\def\psone#1#2{\centerline{
    \epsfxsize #2pc  \epsfbox{#1}
    }}

\begin{document}

\begin{frontmatter}

\title{The Ehrlich-Aberth method for  palindromic matrix polynomials represented in 
the Dickson basis}

\author[dm]{Luca Gemignani}
\ead{gemignan@dm.unipi.it}

\author[dm]{Vanni Noferini\corref{cor1}}
\ead{noferini@mail.dm.unipi.it}

\cortext[cor1]{Corresponding author.}

\address[dm]{Department of Mathematics, University of Pisa\\
Largo Bruno Pontecorvo 5, 56127 Pisa, Italy}

\begin{abstract}
An algorithm based on the Ehrlich-Aberth root-finding method is presented for the
computation of the eigenvalues of a 
T-palindromic matrix polynomial. A structured linearization of the polynomial
represented in the Dickson basis  is introduced in order to
exploit the symmetry of the roots  by halving  the total number of the required  
approximations. 
The rank structure properties of the linearization  allow the design  of a fast  and
numerically robust implementation of the 
root-finding iteration.  Numerical experiments that confirm the effectiveness and
the robustness of the approach are provided. 


\end{abstract}


\end{frontmatter}

\section{Introduction}

The design of effective numerical methods for solving structured generalized
eigenvalue problems has recently 
attracted a great deal of attention. Palindromic matrix polynomials  arise in many
applications \cite{MMMM}. 
An $n\times n$ matrix polynomial of degree $k$ $P(z)=\sum_{i=0}^k A_i z^i$, $A_k\neq
0$, 
$A_i\in \mathbb C^{n\times n}$, $0\leq i\leq k$, is said to be T-palindromic if
$A_i^T=A_{k-i}$ for $i=0, \ldots, k$.    
It is well-known \cite{MMMM}, \cite{prev} that the palindromic structure induces certain  spectral symmetries: 
in particular if $\lambda\neq 0$ is an eigenvalue of $P(z)$  then
$1/\lambda$  is also an eigenvalue  of $P(z)$.  Numerical  solution methods 
 are generally asked to preserve  these symmetries. 

The customary approach  for polynomial eigenproblems consists in two steps: First
$P(z)$ is linearized 
into a matrix  pencil $L(z)=z X + Y$, $X, Y\in \mathbb C^{n k\times n k}$, and then
the eigenvalues of 
$L(z)$ are computed by some iterative solver. The usual choice of the  
matrix QZ algorithm applied to a companion linearization \cite{bible} of 
$P(z)$  is implemented in 
the Matlab function {\em polyeig}.  An alternative solver   based on the
Ehrlich-Aberth root finding algorithm is 
proposed in \cite{BGT} for dealing with certain structured linearizations. 
Specifically the method of \cite{BGT} 
is designed to solve   generalized tridiagonal eigenvalue problems but virtually, as
shown below, it can be extended to 
several other rank structures.   A generalization for tridiagonal  quadratic
              eigenvalue problems is  presented in \cite{bor}. 
A similar strategy using Newton's iteration directly
applied to compute the zeros of 
$P(z)$ is pursued in \cite{WG}.

Modified methods  for palindromic eigenproblems  which are able to preserve their 
spectral symmetries 
have been proposed in several papers. The construction  of T-palindromic
linearizations of palindromic eigenproblems
is the subject of \cite{smith} and \cite{odd},   whereas numerical methods based on
matrix iterations  have been devised in 
\cite{jacobi}, \cite{palqr}, \cite{urv} and \cite{hybrid} for computing the
eigenvalues of these linearizations 
by maintaining the palindromic structure throughout the computation.  To date,
however, the authors 
are not aware of any  specific adaptation of the 
root-finding based methods  to palindromic structures.

The contribution of this paper is to fill the gap by developing a root finder
specifically  suited 
for T-palindromic matrix polynomials, with particular emphasis on the case of large degree.  T-palindromic polynomials of large even degree arise as truncation of Fourier series in several applications such as spectral theory, filtering problems, optimal control and multivariate discrete time series prediction \cite{wilson}.

The polynomial root-finding paradigm 
is a flexible, powerful and quite general tool for solving both structured and unstructured polynomial eigenproblems. 
In its basic form  it proceeds in four steps: 
\begin{enumerate}
\item The matrix polynomial is represented in some convenient polynomial basis. 
\item The transformed polynomial is linearized. 
\item The linearization is reduced in the customary  Hessenberg-triangular form.
\item  A root-finding method is applied for approximating the  eigenvalues of the (reduced) pencil. 
\end{enumerate}
This scheme has some degrees of freedom  concerning the choice of the polynomial basis at step 1 
and the choice of the linearization at step 2 which can be used to exploit both structural and root  properties of the 
matrix polynomial.  The complexity heavily depends on the efficiency of the polynomial zero-finding method
applied to the determinant of the pencil. Steps  2 and 3 are optional  but can substantially improve the 
numerical and computational properties of the method.  Some caution should be used  at step 1 since the change of the basis 
could  modify the spectral structure of the matrix polynomial. The key idea we propose for the implementation of step 4 is the use of the Jacobi formula \cite{golberg}. We emphasize that, although in this paper we focus on palindromics and on a version of the method that is able to extract the palindromic spectral structure, this strategy may be used to address the most general case of an unstructured matrix polynomial eigenproblem, for instance by applying it to the companion linearization. An analysis of the application of the method to a generic matrix polynomial will appear elsewhere.

In this paper we consider  the  polynomial root-finding paradigm  for solving T-palindromic eigenproblems. In particular,  we address 
the main theoretical and computational issues arising  at steps 1, 2 and 4 
of the previous scheme  applied  to  T-palindromic matrix polynomials, and also we indicate  briefly how to  carry out the reduction at step 3. 
The proposed approach relies upon  the
representation and 
manipulation of T-palindromic matrix polynomials  in a 
different degree-graded  polynomial basis $\{\phi_j(y)\}$, namely the Dickson basis,
satisfying a 
three-term recurrence relation and defined by 
$\phi_0(y) = 2$, $\phi_1(y) = y$ and  $y \phi_j(y)=\phi_{j+1}(y)+\phi_{j-1}(y)$ for
$j=1, 2, \ldots$.
For the  given T-palindromic polynomial $P(z)$ of degree $k=2 h$ we determine a
novel polynomial 
$M(y)=\sum_{j=0}^{h+1}M_j \phi_j(y)$, $M_j\in \mathbb C^{2n \times 2n}$, $0\leq
j\leq h+1$, $y=z+z^{-1}$, 
 with the property that   if $\lambda$ and $\lambda^{-1}$ are two distinct 
(i.e. $\lambda \neq \pm 1$) finite semi-simple eigenvalues of $P(z)$ with
multiplicity $\ell$, then $\mu=\lambda+\lambda^{-1}$ 
is a semi-simple eigenvalue for $M(y)$ with multiplicity $2 \ell$.  Moreover, we
find that 
\[g(y)=\det(M(y))=[\det(z^{-h} P(z))]^2 =p(y) \cdot p(y), \]
where $p(y)$ is a polynomial of degree $n h $ at most.

Solving the algebraic equation $p(y)=0$ is at the core of our method for
T-palindromic eigenproblems. 
Our computational experience in polynomial root-finding indicates that the
Ehrlich-Aberth  method \cite{aberth}, \cite{ehrlich} 
for the simultaneous 
approximation of polynomial zeros realizes a quite good  balancing between the
quality of convergence and the cost per iteration.  
The main requirements for the effective implementation of the  Ehrlich-Aberth 
method are both a fast, robust and stable 
procedure to evaluate the Newton correction $p(y)/p'(y)$ and a reliable criterion to
stop the iteration. 
Concerning the first issue  it is worth noting that $p(y)/p'(y)=2 g(y)/g'(y)$  and,
therefore, the computation  
immediately reduces  to  evaluating the Newton correction of $g(y)$.  A suitable
structured linearization $L(y)$ of 
$M(y)$ can be obtained following \cite{bases}  which displays  a semiseparable
structure. In this way,    in  view of  the 
celebrated Jacobi Formula\cite{golberg}  $g'(y)/g(y)=\trace((L(y))^{-1}L'(y))$,  the Newton
correction can be evaluated 
by performing a QR factorization of $L(y)$, say $L(y)= Q(y)\cdot R(y)$,  at  low
computational cost and fulfilling the 
desired requirements of  robustness and stability.   Also,  since 
$\parallel (L(y))^{-1}\parallel_2=\parallel (R(y))^{-1}\parallel_2$ for $y\notin
\spec(L(y))$ we obtain  at no additional cost 
a reliable stop condition based on an estimate of the backward error given by Higham
and Higham \cite{hihi}. 

If $k$, the degree of the matrix polynomial, 
is large with respect to $n$, the size of its matrix coefficients, 
our  approach looks appealing since with a smart choice of the starting points it 
 needs $\mathcal{O}(n^4 k + n^3 k^2)$ operations, whereas the QZ method makes use of 
$\mathcal{O}(n^3k^3)$ operations.  The unpleasant factor $n^4$ in our cost estimate  
depends on the block structure of the 
linearization  used in our current implementation and can in principle  be  decreased  by performing the  
preliminary reduction of 
the linearization in Hessenberg-triangular form as stated at step 3 of the basic scheme. 
The reduction can be carried out by a structured method exploiting the semiseparable  structure of the block 
linearization to compute a rank-structured 
 Hessenberg-triangular linearization.  Incorporating the structured method in our implementation  would finally 
 lead to a 
fast method that outperforms the QZ algorithm for large degrees  and is  comparable in cost for small degrees.

The paper is organized as follows.  The  theoretical  properties of the considered
linearizations of T-palindromic matrix polynomials expressed in the 
Dickson basis are  investigated in Section \ref{teo} and  \ref{teo1}.  The
derivation of the  proposed eigenvalue method   for T-palindromic 
eigenproblems is established in Section \ref{main1} and \ref{main2}.  The complete
algorithm  is described in Section \ref{alg}. Numerical experiments are presented  in 
Section \ref{exp} to illustrate the robustness of our implementation  and to 
indicate computational issues and possible improvements of our algorithm  
compared with other existing methods. 
Finally, conclusion
and future work are discussed in Section \ref{end}.

\section{Theoretical preliminaries on polynomial bases linearizations}\label{teo}

This preparatory section recalls  some basic definitions, background facts
and notations used throughout the paper. 

For $j=0,\dots,k$ let $P_j \in \mathbb{C}^{n \times n}$, $P_k\neq 0$,  be constant
matrices and consider the matrix polynomial 
$P=P(\lambda)=\sum_{j=0}^{k} P_j \lambda^j$. The  generalized polynomial 
eigenproblem  (PEP) associated to $P(\lambda)$ is to find 
an eigenvalue  $\lambda_0$ and a corresponding  nonzero eigenvector  $\B x_0$ 
satisfying 
\begin{equation}
P(\lambda_0) \B x_0  = \B 0.
\end{equation}
In this paper, we will always suppose that $P(\lambda)$ is \emph{regular}, i.e. its
determinant does not identically vanish.

A \emph{linearization} of $P(\lambda)$ is defined as a pencil $L(\lambda)=\lambda X
+ Y$, with $X, Y \in \mathbb{C}^{kn \times kn}$, 
such that there exist unimodular  polynomial matrices $E(\lambda)$ and $F(\lambda)$
for which
$$E(\lambda)L(\lambda)F(\lambda)=\left(\begin{array}{cc}
P(\lambda) & 0\\
0 & I_{(k-1)n}
\end{array}\right).$$
Moreover, if one defines the reversal of a matrix polynomial as 
rev$(P):=\lambda^{k} \sum_j P_j \lambda^{-j}$, the linearization is said to be
\emph{strong} whenever rev$(L)=\lambda Y + X$ is a linearization of rev$(P)$.

Following the work of Mackey, Mackey, Mehl and Mehrmann \cite{prev}, in the paper
\cite{sym} Higham, Mackey, 
Mackey and Tisseur study the two (right and left) \emph{ansatz vector linearization
spaces}: having introduced 
the vector $\B \Lambda:=(1,\lambda,\dots,\lambda^{k-1})^T$, these spaces are defined as
follows:
\begin{eqnarray}
\mathcal{\hat{L}}_1 := \{L=\lambda X+ Y : \exists \B v \in \mathbb{C}^k s. t. L\cdot
(\B \Lambda \otimes I_n) = \B v \otimes P\}\\
\mathcal{\hat{L}}_2 := \{L= \lambda X+ Y : \exists \B w \in \mathbb{C}^k s. t. (\B
\Lambda^T \otimes I_n)\cdot L = \B w^T \otimes P\}.
\end{eqnarray}
It is shown  in \cite{prev} that almost every pencil in these spaces is a
linearization, while in 
\cite{sym}  two binary operations on block matrices, called column shifted sum and
row 
shifted sum, are  first introduced and then  used   to characterize the above
defined spaces.

On the other hand, in \cite{bases} Amiraslani, Corless and Lancaster consider
linearizations of a matrix polynomial 
expressed in some polynomial bases different than the usual monomial one. Equation
(7) in \cite{bases} resembles 
closely the defining equation of $\mathcal{\hat{L}}_2$. The authors themselves
stress this analogy, that suggests an 
extension of the results of \cite{sym} to the case of different polynomial bases. 
Let $\{\phi_i\}_{i=0,\dots,k}$ be a basis for the polynomials of degree less than or
equal to $k$. In \cite{bases}  
degree-graded bases that satisfy a three-terms recurrence relation (for instance,
orthogonal polynomials always do so) are  considered:
\begin{equation}\label{gen}
\lambda \phi_j(\lambda) = \alpha_j \phi_{j+1}(\lambda) + \beta_j \phi_j(\lambda) +
\gamma_j \phi_{j-1}(\lambda).
\end{equation}
The $\alpha_j$ are obviously linked to the leading-term coefficients of the
$\phi_j$. Specifically, calling $c_j$ 
such coefficients, one has that $c_j = \alpha_j c_{j+1}$.

We wish to consider the expansion of the polynomial $P(\lambda)$ in  this basis:
\begin{equation}\label{pgen}
P(\lambda)=\sum_{j=0}^k A_j \phi_j(\lambda).
\end{equation}
We introduce the vector 
\[
\B \Phi:=(\phi_0(\lambda),\phi_1(\lambda),\dots,\phi_{k-1}(\lambda))^T.  
\]
By generalizing the linearizations studied in \cite{bases}, for each choice of $\B \Phi$ two new ansatz
vector linearization spaces can be defined:
\begin{eqnarray}
\mathcal{L}_1 := \{L=\lambda X+ Y : \exists \B v \in \mathbb{C}^k s. t. L\cdot (\B
\Phi \otimes I_n) = c_{k-1} \B v \otimes P\};\\
\mathcal{L}_2 := \{L= \lambda X+ Y : \exists \B w \in \mathbb{C}^k s. t. (\B \Phi^T
\otimes I_n)\cdot L = c_{k-1} \B w^T \otimes P\}.
\end{eqnarray}
It is worth noticing that it is not strictly necessary for the new basis to be
degree-graded, nor it is to satisfy a 
three-term recurrence relation. In fact, it is sufficient that
$\{\phi_i\}_{i=0,\dots,k-1}$ are linearly independent and have 
degree less than  or equal to $k-1$, so that  there exists an invertible  basis
change matrix $B$ such that $\B \Phi=B\B \Lambda$. The basis is degree-graded if and only if $B$ is lower triangular.

In the light of the above definitions  it  is immediately  seen  that  
the main results  of \cite{prev}, \cite{sym} remain valid in the case of a more
general polynomial basis. In particular the following result holds. 

\begin{proposition}\label{prop2}
Let $L \in \mathcal{L}_1$ $(\mathcal{L}_2)$. The following properties are equivalent:
\begin{itemize}
\item $L$ is a linearization of $P$
\item $L$ is a strong linearization of $P$
\item $L$ is regular
\end{itemize}
\end{proposition}
\emph{Proof.} It is a corollary of Theorem 4.3 of \cite{prev}. In fact, any $L\in
\mathcal{L}_1$  $({\rm resp.}, \ \mathcal{L}_2)$ 
can be written as $L=c_{k-1} \hat{L}\cdot  (B^{-1} \otimes I_n)$  $({\rm resp.}, \ L=c_{k-1} (B^{-T}
\otimes I_n) \hat{L})$
for some $\hat{L} \in \mathcal{\hat{L}}_1$ $({\rm resp.}, \  \mathcal{\hat{L}}_2)$. 
Therefore, $L$ has each of the three properties above if and only if $\hat{L}$ has
the corresponding property. $\Box$

This proposition guarantees that almost every (more precisely, all but a closed
nowhere dense set of measure zero) pencil in $\mathcal{L}_1$  $(\mathcal{L}_2)$ 
is a 
strong linearization for $P$. For a proof, see Theorem 4.7 of \cite{prev}.
The eigenvectors of $L$ are related to those of $P$. More precisely, $(\lambda, \B
\Phi \otimes \B x)$ is an eigenpair 
 for $L$ if and only if $(\lambda, \B x)$ 
is an eigenpair  for $P$. Moreover, if $L$ is a linearization then every eigenvector
of $L$ is of the form $\B \Phi \otimes \B x$ for some eigenvector $\B x$ of $P$. 
A similar recovery property holds for the left ansatz vector linearizations. These
properties can be simply proved as in Theorems 3.8 and 3.14 of \cite{prev}, 
that demonstrate them for the special case $\B \Phi = \B \Lambda$.

For the numerical treatment of palindromic  generalized eigenproblems a crucial role
is played by the so-called Dickson  basis \cite{dickson} $\{\phi_i\}_{i\geq 0}$
defined  by  
\begin{equation}\label{dikdik}
\begin{cases}
\phi_0(y) = 2\\
\phi_1(y) = y\\
\forall j\geq 1,~~ y \phi_j(y)=\phi_{j+1}(y)+\phi_{j-1}(y).
\end{cases}
\end{equation}
If we  consider the mapping  $y\colon =\lambda + \lambda^{-1}$  (which we will refer to as the Dickson transformation or the Dickson change of variable)
 then 
$\lambda^j + \lambda^{-j} = \phi_j(y)$  for $j=0, 1, \ldots$.
For  $\lambda=e^{i \alpha}$, we obtain that  $\phi_j(y)=2 \cos(j \alpha)$. 
From  \cite{bases}  by  choosing $\B e_k$ as the ansatz vector    we find   a suitable 
strong 
linearization  of $P(\lambda)$ represented  as in \eqref{pgen}:
{\scriptsize{
\begin{equation}\label{complike}
\left( \begin{array}{cccccccc}
I_n & & & & \\
 & I_n & & & \\
& & \ddots & &\\
& & & I_n & \\
& & & & A_k
\end{array}\right) \lambda + \left( \begin{array}{cccccc}
0 & -2 I_n & & & & \\
-I_n & 0 & -I_n & & & \\
& -I_n & 0 & -I_n & & \\
& & \ddots & \ddots & \ddots &\\
& & & -I_n & 0 & -I_n\\
A_0& A_1 & \dots & A_{k-3} & A_{k-2} - A_{k} & A_{k-1}
\end{array}\right).
\end{equation}}}
In the next section we study  the spectral modifications  induced  by the Dickson 
change of variable  that   provide the
basic link between 
  palindromic matrix polynomials  and  
 matrix polynomials  expressed in the Dickson basis.

\section{Preservation of Jordan structure in the Dickson transformation}\label{teo1}
Let us recall that if $\lambda_0$ is an eigenvalue of  $P(\lambda)$ 
then  the set $\{\B x_j\}$, $j=0,\dots,\ell$ is a Jordan chain of length $\ell +1$
if $\B x_0 \neq \B 0$ and the following relations hold \cite{bible}:
\[
\sum_{i=0}^m \frac{P^{(m-i)}(\lambda_0)}{(m-i)!} \B x_i = \B 0, \quad  m = 0, \dots,
\ell,
\]
where $P^{(k)}(\lambda_0)$ denotes the $k$-th derivative of $P(\lambda)$ evaluated
at $\lambda=\lambda_0$. 
The case $m=0$ corresponds to the definition of an eigenvector. The notion of a
Jordan chain can be 
extended to any matrix function $F : \mathbb{C} \rightarrow \mathbb{C}^{n \times n}$
whose determinant 
vanishes at $\lambda_0$, as long as $F(\lambda)$ is analytic in a neighborhood of
$\lambda_0$. 
In particular, the case of Laurent polynomials is important for our investigations.
If the principal 
part of a Laurent polynomial $L(\lambda)$ is a polynomial of degree $k$ in
$1/\lambda$, then $P(\lambda)=\lambda^k L(\lambda)$ 
is a polynomial. The following lemma relates the Jordan chains of the two. The proof
is a straightforward application of 
the product differentiation rule. 

\begin{lemma}\label{laurent}
Let $L(\lambda)$ be a (Laurent) polynomial and $P(\lambda) = \lambda^k L(\lambda)$
for some natural number $k$. 
Then the set $\{\B x_j\}$ is a Jordan chain of length $\ell +1$ for $P(\lambda)$
associated to the eigenvalue 
$\lambda_0 \neq 0$ if and only if $\{\B x_j\}$ is a Jordan chain of length $\ell +1$
for $L(\lambda)$ associated to the same eigenvalue.
\end{lemma}

Roughly speaking, Lemma \ref{laurent} makes us able to switch between regular and
Laurent polynomials without worrying 
about changes in eigenvalues and generalized eigenvectors. Actually, this result can
be slightly generalized with the next 
lemma, which is just an adaptation  of a well-known  result  in \cite{bible} for the
case where the four matrix functions  that 
we are going to consider are polynomials. In order to prove the lemma, we recall \cite{bible}
that a vector polynomial $\B \phi(\lambda)$ 
is called a root polynomial of order $\ell+1$ corresponding to $\lambda_0$ for the matrix polynomial $P(\lambda)$ if the following conditions are satisfied:
\begin{equation}
\begin{cases}
\phi(\lambda_0) \neq \B 0;\\
\lambda_0 \ \mathrm{is} \ \mathrm{a} \ \mathrm{zero} \ \mathrm{of} \ \mathrm{order} \ \ell+1 \ \mathrm{for} \ P(\lambda) \B \phi(\lambda).
\end{cases}
\end{equation}

Obviously a root polynomial of order $\ell+1$ is defined up to an additive term of 
the form $(\lambda-\lambda_0)^{\ell+1} \B v(\lambda)$ for any suitable vector polynomial $\B
v(\lambda)$. It is possible to prove that $\B \phi(\lambda)=\sum_{j=0}^{\ell}
(\lambda-\lambda_0)^j \B \phi_j + (\lambda-\lambda_0)^{\ell+1} \B v(\lambda)$ if and only if $\{\B \phi_j\}$ is a Jordan 
chain of length $\ell +1$ for $P(\lambda)$ at $\lambda=\lambda_0$. When $\lambda_0 \neq 0$, thanks to Lemma 
\ref{laurent} it is possible to extend the concept to Laurent polynomials: if
$L(\lambda)$ is a Laurent polynomial whose singular 
part has degree $k$ as a polynomial in $\lambda^{-1}$, then we say that $\B
\phi(\lambda)$ is a root polynomial for $L(\lambda)$ 
if it is a root polynomial for $\lambda^k L(\lambda)$.

\begin{lemma}\label{analytic}
Let $P_1(\lambda)$, $P_2(\lambda)$ be (Laurent) polynomials and $A(\lambda)$,
$B(\lambda)$ be 
two matrix functions with $P_2(\lambda)=A(\lambda)P_1(\lambda)B(\lambda)$. Suppose
that an open 
neighborhood $\Omega$ of $\lambda_0\neq 0$  exists 
such that all the considered functions are analytic in $\Omega$, and also suppose
that both 
$A(\lambda_0)$ and $B(\lambda_0)$ are invertible. Then $\lambda_0$ is an eigenvalue
for 
$P_1$ if and only if it is an eigenvalue for $P_2$, and $\{\B y_i\}$ is a Jordan
chain of length 
$\ell +1$ for $P_2$ at $\lambda_0$ if and only if $\{\B z_i\}$ is a Jordan chain of
length $\ell +1$ 
for $P_1$ at $\lambda_0$, where $\B z_i = \sum_{j=0}^i \frac{B^{(j)}(\lambda_0)}{j!}
\B y_{i-j}$.
\end{lemma}
\emph{Proof.} If  $P_1(\lambda)$ and $P_2(\lambda)$ are  classical polynomials then
the thesis follows as in the proof of 
Proposition 1.11 in  \cite{bible} after having represented $A(\lambda)$ and $B(\lambda)$ by
their Taylor  series expansions. 
To deal with the Laurent case, let $\alpha$ and $\beta$ be the minimal integers such
that $Q_1(\lambda):=\lambda^{\alpha}P_1(\lambda)$ 
and $Q_2(\lambda):=\lambda^{\beta}P_2(\lambda)$ are classical polynomials. Just
follow the previous proof for 
$Q_2(\lambda)=\lambda^{\beta-\alpha}A(\lambda)Q_1(\lambda)B(\lambda)$ and apply
Lemma \ref{laurent}.$\Box$

We are now in the position to prove  a result for the Dickson change of variable
$y=\lambda+1/\lambda$. 
The following proposition shows that the number of Jordan chains and their length at
some eigenvalue 
$y_0$ (for the sake of brevity, we shall use the expression \emph{Jordan structure
at $y_0$}) 
is related to the Jordan structures at $\lambda_0$ and $\lambda_0^{-1}$. 

\begin{lemma}\label{change}
Let $y(\lambda)=\lambda+\lambda^{-1}$ and let $M(y)$ be a polynomial in $y$, so that 
$N(\lambda):=M(y(\lambda))$ is a Laurent polynomial in $\lambda$. 
Let first $y_0 = \lambda_0 + \lambda_0^{-1}$,  $\lambda_0 \neq \pm 1$,  be 
a finite eigenvalue of $M(y)$. Then the Jordan structure of $M(y)$ at $y_0$ is 
equal to the Jordan structure of $N(\lambda)$ at either $\lambda_0$ or $1/\lambda_0$. 
 If on the contrary $\lambda_0=\pm 1$, then there is a Jordan chain of length 
$\ell$ at $M(\pm 2)$ if and only if there is a Jordan chain of length $2 \ell$ at
$N(\pm 1)$.
\end{lemma}
\emph{Proof.} It is obvious that $y_0\in \mathbb C$ is an 
eigenvalue for $M(y)$ if and only if both $\lambda_0$ and 
$\lambda_0^{-1}$ are eigenvalues of $N(\lambda)$. 
Let $M(y)=E(y)D(y)F(y)$, where $D(y)=\textrm{diag}(d_1(y),\dots,d_n(y))$ is the 
Smith form (\cite{bible},\cite{smithform}) of $M(y)$. Define $\hat{E}(\lambda):=E(y(\lambda)),
\hat{D}(\lambda):=D(y(\lambda)), \hat{F}(\lambda):=
F(y(\lambda))$. If $\alpha,\beta,\gamma$ are such that
$\tilde{N}(\lambda)=\lambda^{\alpha+\beta+\gamma} N(\lambda)$, $\tilde{E}(\lambda)=
\lambda^{\alpha}\hat{E}(\lambda)$,
$\tilde{D}(\lambda)=\lambda^{\beta}\hat{D}(\lambda)$ and 
$\tilde{F}(\lambda)=\lambda^{\gamma}\hat{F}(\lambda)$ are polynomials in $\lambda$,
then 
we have the relation
$\tilde{N}(\lambda)=\tilde{E}(\lambda)\tilde{D}(\lambda)\tilde{F}(\lambda)$; 
however, in general  
$\tilde{D}(\lambda)$ 
is not the Smith form of $\tilde{N}(\lambda)$. Nevertheless, it has the form
$\textrm{diag}(\lambda^{k_1}\tilde{d}_1(\lambda),
\dots,\lambda^{k_n}\tilde{d}_n(\lambda))$ where $k_1 \geq k_2 \geq \dots \geq k_n$
and $\tilde{d}_i(\lambda)=\lambda^{\deg(d_i)}d_i(y(\lambda))$. 
In other words, the $\tilde{d}_i(\lambda)$-s are palindromic polynomials with no
roots at $0$ and such that  
$\tilde{d}_{i}(\lambda)$ divides $\tilde{d}_{i+1}(\lambda)$ for $i=1,\dots,n-1$.
Moreover, $y_0$ is a zero 
of multiplicity $n$ for $d_i(y)$ if and only if both $\lambda_0$ and
$\lambda_0^{-1}$ are zeros of multiplicity $n$ 
for $\tilde{d}_i(\lambda)$. To reduce $\tilde{D}(\lambda)$ into a Smith form, we
proceed by steps working on $2 \times 2$ principal submatrices.

In each step, we consider the submatrix $\left(\begin{smallmatrix} \lambda^{\alpha}
\tilde{d}_i(\lambda) & 0\\0 & 
\lambda^{\beta}\tilde{d}_{j}(\lambda)\end{smallmatrix}\right)$, with $i<j$. If
$\alpha\leq \beta$, then do nothing; 
if $\alpha > \beta$, premultiply the submatrix by $\left(\begin{smallmatrix} 1 &
1\\-b(\lambda) & 1-b(\lambda)\end{smallmatrix}\right)$ 
and postmultiply it by $\left(\begin{smallmatrix} a(\lambda) &
-q(\lambda)\\b(\lambda) & \lambda^{\alpha-\beta}\end{smallmatrix}\right)$, 
where $q(\lambda)=\tilde{d}_{j}(\lambda)/\tilde{d}_i(\lambda)$ while $a(\lambda)$
and $b(\lambda)$ are such that $a(\lambda) \lambda^{\alpha} 
\tilde{d}_i(\lambda) + b(\lambda) \lambda^{\beta} \tilde{d}_{j}(\lambda) =
\lambda^{\beta} \tilde{d}_i(\lambda)$; the existence of two such 
polynomials is guaranteed by Bezout's lemma, since $\lambda^{\beta} \tilde{d}_i(\lambda)$ is
the greatest common divisor of $\lambda^{\alpha} \tilde{d}_i(\lambda)$ 
and $\lambda^{\beta}\tilde{d}_{j}(\lambda)$. It is easy to check that both matrices are
unimodular, and that the result of the matrix multiplications 
is $\left(\begin{smallmatrix} \lambda^{\beta}\tilde{d}_i(\lambda) & 0\\0 &
\lambda^{\alpha} \tilde{d}_{j}(\lambda)\end{smallmatrix}\right)$. 
By subsequent applications of this algorithm we thus conclude that the Smith form of
$\tilde{D}(\lambda)$ is $\hat{D}(\lambda)=\textrm{diag}
(\lambda^{k_n}\tilde{d}_1(\lambda),\dots,\lambda^{k_1}\tilde{d}_n(\lambda))$.

It follows that the $i$th invariant polynomial of $M(y)$ has a root of multiplicity
$n_{i}$ at $y_0$ if and only if the $i$th invariant polynomial 
of $\tilde{D}(\lambda)$ has a root of multiplicity $n_i$ at $\lambda_0 \neq \pm 1$ and a root
of multiplicity $n_i$ at $1/\lambda_0$. From Lemma \ref{analytic}, 
the Jordan structures of $\tilde{N}(\lambda)$ are equal to those of
$\tilde{D}(\lambda)$. The thesis follows from the properties of the Smith form and 
from Lemma \ref{laurent}.

Mutatis mutandis, a similar argument can be used to analyze the case of $\lambda=\pm 1$: notice in fact that $(y \pm 2)^{k}$ is a factor of the $i$th invariant polynomial of $M(y)$ if and only if $(\lambda \pm 1)^{2 k}$ is a factor of the $i$th invariant polynomial of $\tilde{D}(\lambda)$. $\Box$

\section{Application to palindromic polynomials}\label{main1}

We will now specialize  our analysis to the case of a matrix polynomial with  palindromic structure.

\begin{remark}\label{oddtoeven}
\rm
In this section, we will only treat the case of even degree palindromic matrix polynomials. Notice in fact that an odd degree palindromic may always be transformed to
an even degree palindromic, either by squaring the variable 
($\lambda=\mu^2$) or by multiplication by $(\lambda+1)I_n$. Potentially, both actions may introduce problems: 
squaring the variable adds an additional symmetry $\{\mu,-\mu\}$ 
to the spectrum while multiplying by $\lambda+1$ increases by $n$ the multiplicity of $-1$ as an eigenvalue.
 
However, the first issue may be solved, after passing to Laurent form, by the use of the change of variable $z=(\mu+\mu^{-1})^2$. See also Remark \ref{evenodd}.

Regarding the latter issue, since one knows that he is adding n times $-1$ there is 
no need to compute it: $n$ of the $(n+1)k$ starting points of the Ehrlich-Aberth iteration 
shall be set equal to $-2$, and there they remain with no further corrections. 
The  shortcoming is  that the Jordan structure at $\lambda=-1$ changes.
\end{remark}

Let $\tilde{P}(\lambda)=\sum_{j=0}^{2k} \tilde{A}_j \lambda^j$ be 
a polynomial of even degree. By Lemma
\ref{laurent}, switching to the Laurent form is not 
harmful for finite nonzero  eigenvalues and the corresponding (generalized)
eigenvectors; we can therefore consider its Laurent counterpart
\begin{equation}\label{skp0}
P(\lambda):=\sum_{j=-k}^k A_j \lambda^j.
\end{equation}

Three different kinds  of palindromic structure can be defined. We say that the
Laurent polynomial is purely palindromic 
(resp., $\star$-palindromic, $\star\in \{T,H\}$) if the following relations hold
between its matrix coefficients:
\[
\begin{cases}
\text{Purely palindromic:} \ \ A_j = A_{-j};\\
\text{$\star$-palindromic:} \ \ A_j = A^\star_{-j}.\\
\end{cases}
\]
It is well-known that the palindromic structure induces certain symmetries of
eigenvalues and eigenvectors: 
in particular if $\lambda_0$ is an eigenvalue, $\B x$ is a right eigenvector and $\B
z^T$ is a left eigenvector, then, denoting complex conjugation with the operator $(\cdot)^*$
\[
\begin{cases}
\text{if P is purely palindromic},~~ P(\frac{1}{\lambda_0}) \B x = \B 0$, ~~$\B
z^{T} P(\frac{1}{\lambda_0}) = \B 0;\\
\text{if P is $T$-palindromic}, ~~P(\frac{1}{\lambda_0}) \B z = \B 0$, ~~$\B x^{T}
P(\frac{1}{\lambda_0}) = \B 0;\\
\text{if P is $H$-palindromic}, ~~P(\frac{1}{\lambda_0^*}) \B z^* = \B 0$, ~~$\B
x^{H} P(\frac{1}{\lambda_0^*}) = \B 0.\\
\end{cases}
\]
In this paper we are primarily interested in the design of an efficient solver for
$T-$palindromic eigenproblems.  
A numerical method will be presented  in  Subsection \ref{tpalin}.  The  proposed
approach can  however 
be described very easily  with purely palindromic polynomials. Thus we first
consider  this case for the sake of clarity. 

\subsection{Purely palindromic polynomials}\label{tpure}

The most obvious way to deal with this kind of palindromicity is via introduction of
the change of variable 
$y=\lambda+\lambda^{-1}$, in order to halve the degree of the polynomial. More
explicitly, one can define $Q(y):=P(\lambda(y))$; 
clearly, the purely palindromic structure of $P(\lambda)$ guarantees that $Q(y)$ is
itself a polynomial in the new variable $y$. 
The next proposition is a simple application of Lemmas \ref{laurent} and
\ref{change}, and it relates eigenvectors and Jordan chains of the two polynomials:
\begin{proposition}
When $\lambda_0 \pm 1$, the Jordan structure of $Q(y)$ at the eigenvalue
$y_0=\lambda_0+\lambda_0^{-1}$ is equal to the 
Jordan structure of $P(\lambda)$ at either $\lambda_0$ or $\lambda_0^{-1}$. If $\lambda_0=\pm 1$, $Q(y)$ has a Jordan chain of length $\ell$ at $y_0=\pm 2$ if and only if $P(\lambda)$ has a Jordan chain of length $2 \ell$ at $\lambda_0=\pm 1$.
\end{proposition}
In particular, the eigenvectors of $Q(y)$ at $y_0$ are exactly the same of the
eigenvectors of $P(\lambda)$ at $\lambda_0$ 
(or equivalently at $\lambda_0^{-1}$, since they are the same).

Albeit very attractive, from a numerical point of view this trick is not very
suitable as soon as one considers a 
high degree polynomial. In fact, the matrix coefficients of $Q(y)$ need to be
computed as linear combinations of the ones of $P(\lambda)$. 
Since the powers of a binomial are involved, the coefficients of these linear
combinations would exponentially grow with the polynomial degree. 
To circumvent  this difficulty, we shall make use of the Dickson polynomials
\eqref{dikdik}. 
The polynomial  $Q(y)$ is readily expressed in terms of the $\phi_j(y)$s since in 
the Dickson basis   the coefficients are just the old ones and therefore no
computation at all is needed, namely, 
\begin{equation}
Q(y) = \frac{A_0}{2} \phi_0 + \sum_{j=1}^k A_j \phi_j(y).
\end{equation}
The associated linearization \eqref{complike}  has several  computational advantages
with respect to  other  customary  linearizations
of $P(\lambda)$.  Its size is $nk$ versus $2nk$, the  spectral
symmetries are preserved and, moreover,  the linearization  displays a  semiseparable 
structure. More precisely, it is of the form $D_0 + D_1 y$ where $D_1$ is identity
plus low rank while $D_0$ is Hermitian plus low rank. 
This kind of structure is preserved under  the QZ algorithm  and it 
may  be exploited for the design of an efficient and numerically
robust root-finder applied to the algebraic equation 
$\det Q(y)=0$.

\subsection{T-palindromic polynomials}\label{tpalin}
Consider now a T-palindromic polynomial of even degree $2k$. We will suppose once
more that neither $0$ nor $\infty$ 
are eigenvalues, so that we can divide by $\lambda^k$ and consider the Laurent form
$P(\lambda)$, which is a T-palindromic 
Laurent polynomial of degree $k$ both in $\lambda$ and in $\lambda^{-1}$. Since the
symmetry $\lambda \leftrightarrow \lambda^{-1}$ 
is still present in the spectrum, we expect that the Dickson basis may still play a
role. However, unlike the purely palindromic 
case, it is not possible to directly express a T-palindromic polynomial as a
polynomial in the variable $y$. In fact, splitting  
$P(\lambda)$ as the sum of its symmetric part and its skew-symmetric part we obtain
that  
\begin{equation}\label{skp}
P(\lambda)=A_0 + \sum_{j=1}^k \left[\frac{A_j + A_j^T}{2} (\lambda^j + \lambda^{-j})
+ \frac{A_j-A_j^T}{2} (\lambda^j - \lambda^{-j})\right]. 
\end{equation}
 If we  introduce the  new variables $y:=\lambda + \lambda^{-1}$ and $w:=\lambda -
\lambda^{-1}$, then $P(\lambda)$ can be 
expressed as a bivariate polynomial in $w$ and $y$ which is always linear in $w$,
that is,
\[
Q(y,w)=P(\lambda(y,w))=:B(y)+w C(y).
\]
The property follows from  \eqref{skp} by   substituting  
\[
\lambda^j + \lambda^{-j} = \phi_j(y), \ \lambda^j - \lambda^{-j} = w \left( \frac{1+(-1)^{j+1}}{2} +
\sum_{\ell=1}^{\lceil j/2\rceil}\phi_{j-2\ell +1}(y) \right), \quad j\geq 1.
\]
Notice moreover that $B(y)$ is a 
symmetric polynomial (that is to say, every matrix coefficient is symmetric), $C(y)$
is skew-symmetric, and the operation of transposition 
corresponds to changing the sign of $w$, that is, 
\[
Q^T(y,w)=P^T(\lambda(y,w))=B(y)-wC(y).
\]
In principle one may think of treating $Q(y,w)$ with available techniques for the
bivariate eigenvalue problem (see e.g. \cite{2par1} and references therein), 
but actually $y$ and $w$ are not independent. They are related by the trigonometric
dispersion relation $w^2=y^2-4$. 
This suggests that it is possible to obtain a univariate polynomial by doubling the
dimensions of the matrix coefficients.
Let us define
\[
M(y) = \left(\begin{array}{cc}
B(y) & w^2 C(y)\\
C(y) & B(y)
\end{array}\right).
\]
Then $M(y)$ is a polynomial  in $y$  of degree $k+1$  at most.
 Moreover, it has the following property: if $\lambda_0$ and $\lambda_0^{-1}$ are
two distinct 
(i.e. $\lambda_0 \neq \pm 1$) finite semisimple eigenvalues of $P(\lambda)$ with
multiplicity $m$, then $y_0=\lambda_0+\lambda_0^{-1}$ 
is a semisimple eigenvalue for $M(y)$ with multiplicity $2m$. To see this, notice
first that
\[
M(y) = \diag(\sqrt{w}I_n, 1/\sqrt{w}I_n) \left(\begin{array}{cc} B(y) & w
C(y)\\wC(y) & B(y)\end{array}\right) \diag(1/\sqrt{w}I_n, \sqrt{w}I_n)
\]
and
\[
\left(\begin{array}{cc} B(y) & w C(y)\\wC(y) &
B(y)\end{array}\right)=\frac{1}{2}\left(\begin{array}{cc} I_n & -I_n\\I_n &
I_n\end{array}\right)
 \left(\begin{array}{cc}Q(y,w) &
    0\\  
0 & Q^T(y,w)\end{array}\right) \left(\begin{array}{cc} I_n & I_n\\-I_n &
I_n\end{array}\right).
\]
Hence, we find that 
\[
M(y)= E(w) \left(\begin{array}{cc}Q(y,w) &
    0\\  
0 & Q^T(y,w)\end{array}\right) E^{-1}(w), \ E(w):=\left(\begin{smallmatrix}
\sqrt{\frac{w}{2}} & -\sqrt{\frac{w}{2}}\\ \sqrt{\frac{1}{2w}} & 
\sqrt{\frac{1}{2w}}\end{smallmatrix}\right) \otimes I_n.
\]

Since, as long as $E(w)$ is defined (that is to say
 $w \neq 0, \infty$ or $\lambda \neq 0, \pm 1, \infty$), $\det(E(w))=1$ then 
\begin{equation}\label{funeq}
\det(M(y))=[\det(Q(y,w))]^2, \quad \forall  \ (y, w)\in \mathbb C\times \mathbb C. 
\end{equation}
 Therefore,  $\lambda_0$ has algebraic multiplicity $m$ for $P(\lambda)$ if and
only if $y_0$ has algebraic multiplicity $2m$ 
for $M(y)$. This gives the factorization 
\begin{equation}\label{funeq1}
\det(M(y))=p(y) \cdot p(y), 
\end{equation}
for a suitable polynomial $p(y)$ having  the zero  $y_0$ of  multiplicity $m$.  
Concerning eigenvectors, if  $\lambda_0$ is semisimple, then let $\B x_j$ (resp. $\B
z_j$), $j=1,\dots,m$ be 
the eigenvectors for $P(\lambda)$ (resp. $P^T(\lambda)$) 
corresponding to $\lambda_0$: it can be easily checked that $\{(w_0 \B x^T_j,\B
x^T_j)^T,(-w_0 \B z^T_j,\B z^T_j)^T\}$, 
where $w_0=\lambda_0+\lambda_0^{-1}$, are two linearly independent eigenvectors for
$M(y)$ corresponding to 
$y_0$. Thus, geometric multiplicity is also $2m$. Indeed,  something more can be
said in the more general case of Jordan chains.

\begin{proposition}\label{jc}
Let $y_0=\lambda_0+\lambda_0^{-1}$ be an eigenvalue of $M(y)$ so that $\lambda_0$
and $\lambda_0^{-1}$ are
 eigenvalues for $P(\lambda)$. If $\lambda_0 \neq 0, \pm 1, \infty$ then the Jordan
structure of $M(y)$ at $y_0$ is 
equal to the union of the Jordan structures of $P(\lambda)$ at $\lambda_0$ and at
$\lambda_0^{-1}$.
\end{proposition}
\emph{Proof.} Since $P(\lambda)$ is T-palindromic, it is clear that the Jordan
structure of 
\[
R(\lambda):=\left(\begin{array}{cc}P(\lambda) &
    0\\  
0 & P^T(\lambda)\end{array}\right)
\]
at either $\lambda_0$ or $\lambda_0^{-1}$ is the union of the Jordan 
structures of $P(\lambda)$ at $\lambda_0$ and at $\lambda_0^{-1}$. Define \[
N(\lambda):=M(y(\lambda))=E(w(\lambda)) R(\lambda) E^{-1}(w(\lambda)).
\] 
The matrix function $E(w)$, defined in the previous page, is analytic everywhere in the $w$ complex plane 
but on a branch semiline passing through the origin. Since by hypothesis 
$w_0 \neq 0$, the branch cut can be always chosen in such a way that $E(w)$ 
is analytic in a neighborhood of $w_0=\lambda_0 - \lambda_0^{-1}$, 
and thus $E(w(\lambda))$ is analytic in a neighborhood of $\lambda_0$. 
Then we can apply Lemma \ref{analytic} to  conclude that the Jordan structures of 
$M(\lambda)$ and $R(\lambda)$ are the same. Application of Lemma \ref{change}
completes the proof. $\Box$
\begin{remark}\label{ins1}
\rm{Another remarkable  property of $M(y)$ is that its coefficients are all
skew-Hamiltonian, that is to say 
they can be written as $JK$ where $J=\left(\begin{smallmatrix} 0 & -I\\ I
&0\end{smallmatrix}\right)$ and 
$K$ is some skew-symmetric matrix. This link between T-palindromic and
skew-Hamiltonian polynomials is 
interesting because it may shed more light on the relation between several polynomial structures. It is 
known that one can easily transform a palindromic polynomial to an even polynomial
by a Cayley transformation, 
and then to a Hermitian polynomial via a multiplication by $i$ (if one started from
a real polynomial) or to a 
symmetric polynomial by squaring the matrix coefficients. On the other hand,
Hamiltonian polynomials can lead to 
skew-Hamiltonian polynomials by squaring each coefficients, and multiplication by
$J$ sends a skew-Hamiltonian 
polynomial to a skew-symmetric polynomial. The Dickson change of variable, followed 
by doubling the dimension, is able to
map T-palindromic polynomials of even degree
 to a subset of skew-Hamiltonian polynomials. Unlike some of the other mentioned
maps, this is not a bijection between
 two classes of structured polynomials, because what is obtained is actually a
subset of skew-Hamiltonian polynomials. In fact, since
 the north-west and south-east coefficients of $M(y)$ are the coefficients of $B(y)$
they must be symmetric and there 
is a relation between the north-east and south-west coefficients of $M(y)$. 
However, a deeper investigation on this subject is needed in the future.}
\end{remark}
\begin{remark}\label{evenodd}
\rm{Notice that a similar technique can be applied to \emph{even/odd} matrix polynomials, that is polynomials whose coefficients alternate between symmetric and skew-symmetric matrices. In this case, on can apply the transformation $z=\lambda^2$ and use algebraic manipulations, akin to the ones described for the T-palindromic case, in order to build a new polynomial in $z$ with double dimensions.

In the case of an odd-degree T-palindromic polynomial, the substitutions $\lambda=\mu^2$ and $y=\mu+\mu^{-1}$ lead to an $M(y)$ such that $\left(\begin{smallmatrix} 0 & I\\ I
&0\end{smallmatrix}\right) \cdot M(y)$ is odd. Therefore, one may apply $z=y^2$ and build a third polynomial in order to extract the additional structure $\{\mu,-\mu\}$.}
\end{remark}
 Equation \eqref{funeq}  and \eqref{funeq1}  enable the computation of the
eigenvalues of $P(\lambda)$ to be reduced to 
solving  algebraic equations. From  Proposition \ref{jc} it follows that 
possible discrepancies in the Jordan structures can be expected for $y_0=\pm 2$  and
$y_0=\infty$ corresponding to $\lambda_0=\pm 1$ and 
$\lambda_0=0, \infty$, respectively.  

When $\lambda_0=\pm 1$ not only the proof we
gave is not valid (because, 
since $w_0=0$ is a branch point, there is no neighborhood of analyticity of the
matrix function $E$), 
but in fact the proposition itself does not hold. As a counterexample, let $a \neq \pm \frac{i}{\sqrt{2}}$ and consider the
polynomial 
\[
P(\lambda)=\left(\begin{array}{cc} \lambda -2 + \lambda^{-1} & a \lambda - a
\lambda^{-1}\\ 
-a \lambda + a \lambda^{-1} & \lambda + \lambda^{-1}\end{array}\right).
\]
 We have that 
$\{\left(1, 0\right)^T,\left(0, a \right)^T\}$ is a 
Jordan chain for $P(\lambda)$ at $\lambda=1$. The corresponding $M(y)$ is 
\[
M(y)=\left(\begin{array}{cccc} y-2 & 0 
& 0 & ay^2-4a\\0 & y & 4a-ay^2 & 0\\ 0 & a & y-2 & 0\\ -a & 0 & 0 &
y\end{array}\right),
\]
which has a 
semisimple eigenvalue at $y=2$ with the corresponding eigenvectors $\left(0, 0, 1,
0\right)^T$ 
and $\left(2,0,0,a\right)^T$.

If the leading coefficient of $P(\lambda)$ is not symmetric, then $M(y)$ has $2n$ extra
infinite eigenvalues, 
where $n$ is the dimension of the matrix coefficients of $P(\lambda)$. These
eigenvalues are 
defective since their geometric multiplicity is only $n+\dim \ker C_{k-1}$, where
$C_{k-1}$ is the leading 
coefficient of $C(y)$.

For the numerical  approximation of the roots of  $p(y)$  we  can exploit again the
properties of the Dickson basis
to compute the matrix coefficients of $M(y)=\sum_{j=0}^{k+1}M_j \phi_j(y)$.  The
code below computes the matrices 
$M_j\in \mathbb C^{2n\times 2n}$, $0\leq j\leq k-1$,   given in input the
coefficients $A_j$ of $P(\lambda)$, $0\leq j\leq k$, 
defined as in \eqref{skp0}.

\begin{code}
{\bf function Dickson$\_$transform}\\
{\bf Input:} $A_0, \ldots, A_k \in \mathbb C^{n\times n}$\\
{\bf Output:} $M_0, \ldots, M_{k+1} \in \mathbb C^{2n\times 2n}$\\
\>$B_0=A_0/2$; $\hat C_0=0_n$; \\
\> for $j=1,\dots,k$ \\
\> \> $B_j=(A_j + A_j^T)/2$; $\hat C_j=(A_j - A_j^T)/2$;\\
\> end\\
\>  $S_0=0_n$,  $S_1=0_n$; \\
\> for $j=k, \ldots ,1$ \\
\> \>  $S_{\mbox{mod}(j, 2)}=S_{\mbox{mod}(j, 2)}+\hat C_{j}$\\
\> \> $C_{j-1}=S_{\mbox{mod}(j, 2)}$; \\
\> end\\
\> $C_0=C_0/2$; $C_k=C_{k+1}=0_n$;  $\tilde C_0=C_2$ \\
\> $\tilde C_1=C_1+C_3$; $\tilde C_2=2 C_1+ C_4$; \\
\> for $j=4,\dots,k$ \\
\> \> $\tilde C_{j-1}=C_{j-3}+C_{j+1}$; \\
\> end\\
\> $\tilde C_k=C_{k-2}$;  $\tilde C_{k+1}=C_{k-1}$; \\
\> for $j=1:k+2$ \\
\> \> $\tilde C_{j-1}=\tilde C_{j-1} - 2 C_{j-1}$; \\
\> \> $M_{j-1}=[B_{j-1}, \tilde C_{j-1}; C_{j-1}, B_{j-1}]$; \\
\> end
\end{code}

\begin{remark}
\rm
The coefficients of $C(y)$ are linear combinations of $A_j - A_j^T$. As can be seen by the above algorithm, the coefficents of such combinations expressed in the Dickson basis remain bounded, the upper bound being $1/2$. An analogous result, with upper bound $1$, holds for $w^2 C(y)$. This is in contrast with the exponential growth that would have been seen in the purely palindromic case if one had directly applied the Dickson transformation without the use of the Dickson basis.
\end{remark}
The arithmetic cost is $\mathcal{O}(n^2k)$ operations. Once the coefficients $M_j$ are
determined,  a linearization of 
$M(y)$ of the form  \eqref{complike} can  be constructed. The properties of this
linearization are investigated in the 
next section in order to devise a fast and numerically robust method to  evaluate
the Newton correction  of $p(y)$ defined by \eqref{funeq1}.

\section{Computing the Newton correction}\label{main2}
Our aim in this section is to derive a fast, robust and stable  method for computing
the 
Newton correction $p(y)/p'(y)=2  \det(M(y))/(\det(M(y)))'$, where $p(y)$ and $M(y)$ are related by 
\eqref{funeq1},  given a structured
linearization 
$L(y)=y E + F$, with $E, F \in \mathbb{C}^{2n(k+1) \times 2n(k+1)}$,  of 
$M(y)$ of the form \eqref{complike}, namely,
\[
E=\left( \begin{array}{cccccccc}
I_{2n} & & & & \\
 & I_{2n} & & & \\
& & \ddots & &\\
& & & I_{2n} & \\
& & & & M_{k+1}
\end{array}\right) 
\]
and 
\[
F=\left( \begin{array}{cccccc}
0 & -2 I_{2n} & & & & \\
-I_{2n} & 0 & -I_{2n} & & & \\
& -I_{2n} & 0 & -I_{2n} & & \\
& & \ddots & \ddots & \ddots &\\
& & & -I_{2n} & 0 & -I_{2n}\\
M_0& M_1 & \dots & M_{k-2} & M_{k-1} - M_{k+1} & M_{k}
\end{array}\right).
\]
Our approach relies upon the celebrated {\em Formula of Jacobi}\cite{golberg}
\[
(\det(L(y)))'=\det(L(y)) \trace( L^{-1}(y)L'(y))=\det(L(y)) \trace( L^{-1}(y) E)
\]
which reduces the evaluation  of  $\det(M(y))/(\det(M(y)))'= \det(L(y))/(\det(L(y)))'$
to computing the trace of  $L^{-1}(y)\cdot  E$.  In the sequel we   describe a
method for finding the block 
entries and, a fortiori, the trace  of  the inverse of $L(y)$  from the LQ 
factorization of the  matrix.
Then we slightly modify  the computation to take into account the  contribution due
to the matrix $E$. It will be clear from what follows that this method is general and can be applied, with only trivial modifications, to any kind of unstructured matrix polynomial, simply by considering for instance the standard companion linearization instead of \eqref{complike}.

We denote as  $\mathcal G(\theta, \psi)$ the $2\times 2$ unitary Givens rotation
given by 
\[
\mathcal G(\theta, \psi)=\left(\begin{array}{cc} \theta & \psi \\ -\bar \psi & \bar
\theta\end{array}\right), \quad |\theta|^2 + |\psi|^2=1.
\]
Let $L(y)=\tilde{L}\cdot Q$  be the  (block) LQ factorization of $L(y)$ obtained by means of
Givens rotations so that 
\begin{equation}\label{schur}
\begin{array}{ll}
L(y) \mathcal G_1 \cdot \mathcal G_2 \cdots \mathcal   G_k = \tilde{L}, \quad  Q^H =\mathcal
G_1 
\cdot \mathcal G_2 \cdots \mathcal   G_k, \\
\mathcal G_j=I_{2n(j-1)}\oplus(\mathcal G(\theta_j, \psi_j)\otimes I_{2n})\oplus
I_{2n(k-j)}.
\end{array}
\end{equation}
It can be easily checked that the lower triangular factor $\tilde{L}$  has the following
structure 
\[
\tilde{L}=\left( \begin{array}{cccccc}
\alpha_1 I_{2n} &  & & & & \\
\beta_1 I_{2n} & \alpha_2 I_{2n} & & & & \\
\gamma_1 I_{2n}& \beta_2 I_{2n} & \ddots  &  & & \\
&\ddots & \ddots & \ddots &  &\\
& & \gamma_{k-2}I_{2n}& \beta_{k-1}I_{2n} & \alpha_kI_{2n}& \\
\hat M_0& \hat M_1 & \dots & \hat M_{k-2} & \hat M_{k-1} & \hat M_{k}
\end{array}\right), 
\]
where $\alpha_j\neq 0$, $1\leq j\leq k$. 
If $\hat M_k$ and, therefore, $L(y)$ is invertible then  the  LQ factorization can
be used to  find a condensed representation of the inverse of 
$L(y)$.   Observe that  $L^{-1}(y)= Q^H  \cdot \tilde{L}^{-1}$.  In order  to take into
account the occurrence of the matrix $E$ in the Jacobi formula let us  introduce 
the matrix $\tilde M_{k+1}=\hat M_k^{-1} \cdot M_{k+1}$.  Then we have the following 
\begin{proposition}\label{inv} 
There exist matrices $\tilde M_1, \ldots, \tilde M_k \in \mathbb C^{2n\times 2n}$
such that 
\[
L^{-1}(y) E=\left( \begin{array}{cccc} \tilde M_1 & \psi_1\tilde M_2  & \ldots &
\psi_1\cdots \psi_k \tilde M_{k+1} \\ & \bar \theta_1  \tilde M_2 \\ &  & \ddots \\
& &  & \bar \theta_k  \tilde M_{k+1}
\end{array}\right), 
\]
where the blank entries are not specified. 
\end{proposition}
\emph{Proof.} The proof  basically  follows by applying the (block) Schur
decomposition 
\eqref{schur} of $Q^H$  to  the block lower triangular factor $\tilde{L}^{-1} E$. 
To show it more formally we can proceed by induction.  Let us assume that the 
the $j-$th block row of $ \mathcal G_{j} \cdots \mathcal   G_k \tilde{L}^{-1} E$ 
can be represented as 
\[
\left[\begin{array}{cccccccc}\star  & \ldots & \star & \tilde M_j & \psi_{j}\tilde
M_{j+1}&\ldots &
 \psi_{j}\cdots \psi_k\tilde M_{k+1}\end{array}\right], 
\]
where $\tilde M_{j}$ is the diagonal entry and the value of the 
entries in the strictly lower triangular part -- denoted by $\star$-- 
 is  not essential. 
Then, by applying $ \mathcal G_{j-1}$ on the left of the  matrix we find that the
$(j-1)-th$ block row 
looks like 
\[
\left[\begin{array}{cccccccc}\star  & \ldots &\star  & \tilde M_{j-1} &
\psi_{j-1}\tilde M_{j}&\ldots &
 \psi_{j-1}\cdots \psi_k\tilde M_{k+1}\end{array}\right], 
\]
whereas the diagonal entry in position $j$ becomes $\bar \theta_{j-1} \tilde M_j$. 
 $\Box$

This result says that the  block diagonal entries of $L^{-1}(y)$ can be  determined
from  the entries in its first (block) row.  The computation of this row is 
equivalent to the solution of the linear system 
\[
\left(I_{2n}, 0_{2n}, \ldots, 0_{2n}\right)= \left( X_1, \ldots, X_{k+1}\right)
\cdot L(y)
\]
or, equivalently, 
\[
\left(I_{2n}, 0_{2n}, \ldots, 0_{2n}\right) \cdot Q^H = \left(X_1, \ldots,
X_{k+1}\right) \cdot \tilde{L}.
\]
In the view of the structure of $Q^H$ this reduces  to
\[
\left( \theta_1, \theta_2 \psi_1, \ldots, \theta_k \prod_{j=1}^{k-1}\psi_j, 
\prod_{j=1}^{k}\psi_j\right)\otimes I_{2n}=\left(X_1, \ldots, X_{k+1}\right) \cdot \tilde{L}.
\]
Let $D\in \mathbb C^{2n(k+1)\times 2n(k+1)} $ be a block diagonal matrix defined by 
\[
D=\diag(1, \psi_1, \ldots,  \prod_{j=1}^{k-1}\psi_j, 
\prod_{j=1}^{k}\psi_j) \otimes I_{2n}.
\]
Using the matrix $D$ to balance the coefficient matrix yields 
 \[
\left( \theta_1, \theta_2, \ldots, \theta_k, 
1\right)\otimes I_{2n}=\left(X_1, \ldots, X_{k+1}\right)D^{-1} \cdot D \cdot \tilde{L}
\cdot D^{-1}.
\]
Observe that 
\[
 \left(\hat X_1, \ldots, \hat X_{k+1}\right)=\left(X_1, \ldots, X_{k+1}\right)D^{-1}= 
\left(\tilde  M_1, \ldots, \tilde M_{k-1}, \hat M_k^{-1}\right), 
\]
and, therefore, the solution of 
\[
\left( \theta_1, \theta_2, \ldots, \theta_k, 
1\right)\otimes I_{2n}= \left(\hat X_1, \ldots, \hat X_{k+1}\right) \hat L, \quad D
\tilde{L} D^{-1} =\hat L, 
\]
gives the desired unknown matrices  $\tilde M_1, \ldots, \tilde M_k$. 
To achieve some computational savings we rewrite  the system as
\[
\left( \theta_1, \theta_2, \ldots, \theta_k, 
1\right)\otimes \hat M_k=\left(\tilde  X_1, \ldots, \tilde X_{k+1}\right) \hat L
\]
and thus   we arrive at the following relation 
\[
\det(M(y))'/\det(M(y))=\trace(\hat M_k^{-1}( \tilde  X_1+ \bar \theta_1 \tilde  X_2
+\ldots +\bar \theta_{k-1} \tilde 
 X_k  + \bar \theta_k   M_{k+1})), 
\]
which is used to compute the reciprocal of the Newton correction.  The function
${\bf trace}$  below implements our 
resulting algorithm at the cost of $\mathcal{O}(n^2  k + n^3 )$ operations.

\begin{code}
{\bf function trace}\\
{\bf Input:} $M_0, \ldots, M_{k+1} \in \mathbb C^{2n\times 2n}$, $\lambda\in \mathbb
C$, $(\det(M(\lambda))\neq 0)$\\
{\bf Output:} the value of $\eta=p'(\lambda)/p(\lambda)$ \\
\>$M_{k-1}=M_{k-1}-M_{k+1}$; $M_k=M_k +\lambda M_{k+1}$; \\
\>$\B \alpha=\lambda \  \ones(1,k+1)$;  $\B \beta=-\ones(1,k+1)$;  \\
\>$\B \gamma=\zeros(1,k)$;  $\B \chi=-\ones(1,k+1)$; $\chi_1=-2$;\\
\> for $j=1,\dots,k$ \\
\> \> $\B v=[\alpha_j; \chi_j]$; $\mathcal G^T=\planerot(\B v)$; $\B q(j,:)=\mathcal
G(1,:)$; $c_j=\B q(j,1)$; \\
\> \> $\alpha_j=\alpha_j \mathcal G_{1,1} +\chi_j \mathcal G_{2,1}$;  $\tilde
\beta=\beta_j \mathcal G_{1,1} +\alpha_{j+1} \mathcal G_{2,1}$; \\
\> \> $\alpha_{j+1}=\beta_j \mathcal G_{1,2} + \alpha_{j+1}\mathcal G_{2,2}$;
$\beta_j=\tilde \beta$; $\gamma_j=\beta_{j+1} \mathcal G_{2,1}$;
$\beta_{j+1}=\beta_{j+1} 
\mathcal G_{2,2}$; 
\\
\> \> $\tilde M=\mathcal G_{1,1} M_{j-1} + \mathcal G_{2,1} M_j$; $M_j=\mathcal
G_{1,2} M_{j-1} + \mathcal G_{2,2} M_j$; $M_{j-1}=\tilde M$;\\
\> end\\
\> for $j=1,\dots,k-1$\\
\> \> $\beta_j=\beta_j  \B q(j, 2)$; \\
\> end\\
\> for $j=1,\dots,k-2$\\
\> \> $\gamma _j=\gamma_j  \B q(j, 2)  \B q(j+1,2)$; \\
\> end\\
\> $s=1$;\\
\>  for $j=k,\ldots, 1$\\
\> \> $s=s   \B q(j,2)$; $M_{j-1}= s  M_{j-1}$; \\
\> end\\
\> $\tilde X_k =(c_k  M_k  -M_{k-1})/\alpha_k$;  
$\tilde X_{k-1} =(c_{k-1}  M_k  -M_{k-2} -\beta_{k-1} \tilde X_k)/\alpha_{k-1}$;\\
 \> for $j=k-2, \ldots, 1$\\
\> \> $\tilde X_{j} =(c_{j}  M_k  -M_{j-1} -\beta_{j} \tilde X_{j+1} -\gamma_j
\tilde X_{j+2})/\alpha_{j}$;\\
\> end\\
\> $\tilde M=\tilde X_1$;\\
\> for $j=1, \ldots, k-1$\\
\> \> $\tilde M=\tilde M +\bar c_j \tilde X_{j+1}$; \\
\> end \\
\> $\tilde M=\tilde M +\bar c_k M_{k+1}$;  $\tilde M=M_k \backslash \tilde M$;
$\eta=\trace(\tilde M)$; 
\end{code}

\section{The  Ehrlich-Aberth algorithm for T-palindromic eigenproblems}\label{alg}
A  simple tool for the simultaneous approximation of all the eigenvalues of a
polynomial is the Ehrlich-Aberth  method.  Bini and Fiorentino 
\cite{BF} showed that a careful implementation of the method
yields an efficient and robust  polynomial root finder.  The software package 
MPSolve documented in \cite{BF}  is designed to successfully compute approximations 
of polynomial zeros at any specified accuracy using  a multi-precision arithmetic environment. 

A root finder for T-palindromic eigenproblems can be based on the Ehrlich-Aberth
method 
applied  for the solution of the algebraic equation  $p(y)=0$, where 
$p(y)$ is related with $M(y)$  by \eqref{funeq1} and $M(y)$ is generated by  the
function 
{\bf  Dickson$\_$transform} applied to the input  coefficients $A_j\in \mathbb
C^{n\times n}$ of 
the  T-palindromic matrix polynomial $P(\lambda)$  of degree $2k$  given  as in
\eqref{skp}. 
The method simultaneously approximates all the zeros of  the polynomial $p(y)$: given
a vector 
$\B z^{(0)}\in \mathbb C^{N}$, $N=n k$,  of 
initial approximations to the zeros of $p(y)$, the Ehrlich-Aberth iteration 
generates a sequence 
$\{\B z^{(k)}\}$, $k\geq 0$, which locally converges to the $N-$tuple of the roots
of $p(y)$, 
according to the equation 
\[
z^{(k+1)}_j= z^{(k)}_j-\frac{p(z^{(k)}_j)/p'(z^{(k)}_j)}{1-\frac{p(z^{(k)}_j)}{
p'(z^{(k)}_j)}\sum_{\ell=1, \ell\neq j}^N \frac{1}{ z^{(k)}_j- z^{(k)}_{\ell}}},
\quad 1\leq j\leq N.
\]
The convergence is superlinear for simple roots and linear for multiple roots. 
In practice, the Ehrlich-Aberth method exhibits quite good global convergence properties, even though
no theoretical results are known in this regard. The main requirements for an efficient implementation
of the method are: 
\begin{enumerate} 
\item  a rule for choosing the initial approximations; 
\item a fast, numerically robust and stable method to compute the Newton correction
$p(z)/p'(z)$; 
\item a reliable stopping criterion.
\end{enumerate}

Concerning the first issue    it is commonly advocated\cite{avvocato}
 that  for scalar polynomials 
the  convergence  benefits  from  the choice of equally spaced points lying on some 
circles around the origin in the complex plane. In the case of matrix polynomials where the eigenvalues 
are often   widely varying in magnitude  this choice can not be optimal. 
A better strategy  using the initial guesses lying on certain ellipses around the origin in the complex plane is 
employed in our method.   The second task can be accomplished by 
means of the   function {\bf  trace} in the previous section.  With respect to the
third issue, it is worth 
observing that  the QL-based method pursued for the trace computation  also provides
an estimate on the backward error
for the generalized eigenvalue problem.  From a result in \cite{hihi} it follows
that if $\tilde y$ is not an eigenvalue of $L(y)$ then 
\[
\eta(\tilde y)=1/(\parallel (\tilde y E + F)^{-1} \parallel_2 (1+ | \tilde y|)) 
\]
gives an  appropriate measure of the backward error  for the approximate eigenvalue
$\tilde y$. 
Since  for $\tilde y E + F =\tilde{L} \cdot Q$ we have  that 
\[
\parallel (\tilde y E + F)^{-1} \parallel_2 =\parallel L^{-1} \parallel_2 \geq
\parallel \hat M_k^{-1}\parallel_2\geq (\sqrt{2n})^{-1} 
\parallel \hat M_k^{-1}\parallel_\infty.
\]
In our implementation we consider the 
quantity 
\[
\hat \eta(\tilde y)= \sqrt{2n}/(\parallel  \hat M_k^{-1} \parallel_\infty (1+ |
\tilde y|))
\]
as an error measure. If $\hat \eta(\tilde y)$ is smaller than a fixed tolerance then
$\tilde y$ is taken as an approximate eigenvalue and the 
corresponding iteration is stopped. 
The resulting Ehrlich-Aberth algorithm for approximating finite
eigenvalues of $M(y)$ and hence obtaining the corresponding eigenvalues of $P(\lambda)$  is described below. 
In the next section we present results of numerical experiment assessing the
robustness of the proposed approach. 

\begin{code}
{\bf function palindromic$\_$aberth$\_$zeros}\\
{\bf Input:} $A_0, \ldots, A_k \in \mathbb C^{n\times n}$,  $tol \in \mathbb R$, $maxit\in \mathbb N$,  initial guesses $z_1, \ldots, z_N$ \\
{\bf Output:} approximations $\zeta_1, \ldots, \zeta_{2N}$, $N=nk$, 
 of the zeros of  $P(\lambda)=\sum_{i=-k}^k A_i \lambda^i$ \\
\> $[M_0, \ldots, M_{k+1} ]=${\bf {Dickson$\_$transform}}$(A_0, \ldots, A_k)$\\
\> $N=nk$; $\B c=\ones(N,1)$; \\
\> $nn=0$;\\
\> for $i=1, \ldots, maxit$\\
\> \> \> for $j=1, \ldots, N$\\
\> \> \> \> if $(c(j))$\\
\> \> \> \> \> $z=${\bf {trace}}$(M_0, \ldots, M_{k+1}, z_j)$;  $z=2/z$; \\
\> \> \> \> \> $h=\pippo(1./(z(1\colon j-1)-z(j)))$;  \\
\> \> \> \> \>$h=h+\pippo(1./(z(j+1\colon N)-z(j)))$;\\ 
\> \> \> \> \> $h= z/(1+h\ z)$; $z_j=z_j-h$; \\
\> \> \> \> \> if $(\hat \eta(z_j)\leq tol \  {\emph or} \  |h|\leq tol |z_j|)$\\
\> \> \> \> \> $c(j)=0$; $nn=nn+1$; \\
\> \> \> \> end;\\
\> \> \> end \\
\> \> end\\
\> \> if $(nn=N)$ \\
\> \> \> break\\
\> \> end\\
\> end\\
\> for $j=1, \ldots, N$\\
\> \> $\B r=\roots([1, -z_j, 1])$;\\
\> \>  $\zeta_{2j-1}=r_1$; $\zeta_{2j}=r_2$;\\
\> end
\end{code}
\medskip

The total cost of the algorithm is therefore $\mathcal{O}(t(n^2  k + n^3))$ operations, 
where $t$ is the total number of times that the function ${\bf trace}$ is called. 
Numerical experiments presented in the next section show that $t$ heavily depends 
on the choice of the starting points. With a smart choice, $t$ is of order 
$\mathcal{O}(nk)$, which gives a total computational cost of $\mathcal{O}(n^4 k + n^3 k^2)$. 
Since the cost of our method grows as $n^4$ but is only quadratic in $k$, 
where customary QZ-like methods use $\mathcal{O}(n^3 k^3)$ operations, an Ehrlich-Aberth approach 
looks particularly suitable when the matrix polynomial has a high degree and small coefficients so that $k^2/n$ is large.

It is worth noticing that the case of large $n$ can still be treated by means of an Ehrlich-Aberth method in $\mathcal{O}(n^3k^3)$ operations. 
  The basic  observation is that  the factor $n^4$ comes from 
the block structure of the linearization involved in the computation of the trace.   A reduction of the cost  can therefore 
 be achieved by a different strategy  where  the linearization is initially converted    into (scalar) triangular-Hessenberg form: say, $N(y)=R y + H$ where $R$ 
is (scalar) triangular and $H$ is (scalar) Hessenberg. The task can virtually be performed by any extension of the fast structured methods 
for the  Hessenberg reduction proposed in \cite{DB,eid}.  These methods  preserve the rank structure  which can therefore be exploited also  
in the triangular-Hessenberg linearization. Once the matrices $R$  and $H$  have been determined then the computation of $\mathrm{tr}(N(y))$ 
can be performed by the following algorithm  which has a cost of $\mathcal{O}(n^2 k^2)=\mathcal O(N^2)$ operations:
\begin{itemize}
\item Perform a $R Q$  decomposition of the Hessenberg matrix $N(y)$, obtaining a unitary matrix $Q$ represented as product of $\mathcal O(N)$ Givens transformations
(Schur decomposition)  and a triangular matrix $U$.
\item Compute the last row of $N(y)^{-1} R$ by solving $\B w^T N(y)=\B e_{N}^T $ and then computing $\B w^T\colon=\B w^T R$.  
\item Recover the diagonal entries of $N(y)^{-1} R$ from the entries of $\B w$ and the  elements of the Schur decomposition of $Q$. 
\end{itemize}

This alternative road leads to an algorithm of total cost $\mathcal{O}(n^3 k^3)$ operations. 
An efficient implementation exploiting the rank structures 
of the matrices involved will be presented elsewhere.

\section{Numerical Experiments}\label{exp}
The function {\bf palindromic$\_$aberth$\_$zeros} for computing the roots of a T-palindromic matrix
polynomial $P(\lambda)=\sum_{j=-k}^k A_j \lambda^j$, 
given its coefficients 
$A_{-j}=A_j^T\in \mathbb C^{n\times n}$, $0\leq j\leq k$,  has been implemented in  
Matlab\footnote{Matlab is a registered trademark of The Mathworks, Inc..} 
and  then    
used for the 
computation of the  zeros of polynomials of both small and high  degree.  
The tolerance is fixed at $tol=1.e-13$ and  for the maximum number of iterations we
set 
$maxit=2nk$.

Extensive numerical experiments have been performed to illustrate some  basic issues 
 concerned with  the efficiency and  the accuracy  of a practical implementation 
of our method.  

\subsection{Efficiency of root-finding}
An accurate and efficient root-finder is essential to the success 
of our algorithm. In practice, the cost of each iteration is 
strongly dependent on the amount of early convergence (for the sake of brevity, in the following we will refer to this phenomenon using the word \emph{deflation}) occurring for a given problem.  In other words,  
a critical point to assess the efficiency of the novel method is the evaluation of the total number $t$ of calls 
of the function trace, and of its dependence on the total number $N:=nk$ of the eigenvalues. 
When the Ehrlich-Aberth method is used to approximate scalar polynomials roots, experiments 
show that $t$ depends on the choice of the starting points. 
If there is not any a priori knowledge about the location of the roots, 
empirical evidence \cite{BF} shows that choosing starting points distributed on some circles around the origin  
leads to  acceptable  performances and/or quite 
regular convergence patterns.  

The class 
 H$_{n,k}$ of   T-palindromic polynomials have been used to  verify  if these properties 
still hold in the matrix case.  
 The polynomials  are  constructed according to the following rules:  
\[
\begin{array}{ll}
{\rm H}_{n, k}=\sum_{j=-k}^k A_j \lambda^j, \quad A_j\in \mathbb R^{n\times n}, \\
A_0=0_n; \  A_j=I_n + \B e_n \B e_1^T,  \ A_{-j}=A_j^T, \  1\leq j\leq k.
\end{array}
\]
From  
\[
h(\lambda)=\sum_{j=1}^k \lambda^j + \sum_{j=1}^k
\lambda^{-j}=\frac{\lambda^{k}-1}{\lambda-1}\frac{\lambda^{k+1}+1}{\lambda^k}, 
\]
we find that most of the eigenvalues lie on the unit circle and for $k$ even
$\lambda=-1$ is a double root of $h(\lambda)$.

Figure 1  describes  the convergence history for our  root finder applied to
$H_{5, 20}$ with  starting values equally spaced on 
the circle centered in the origin with radius 4.  The curves represented  are
generated by plotting the sequences $\{z_j^{(k)}\}$, $1\leq k \leq maxit$, 
for $j=1, \ldots, N$.  The convergence is quite regular and very similar to that
exhibited in the scalar polynomial case \cite{BF}   and  theoretically predicted for 
simultaneous iterations based on Newton-like methods \cite{Hub}.

\begin{center}
\psone{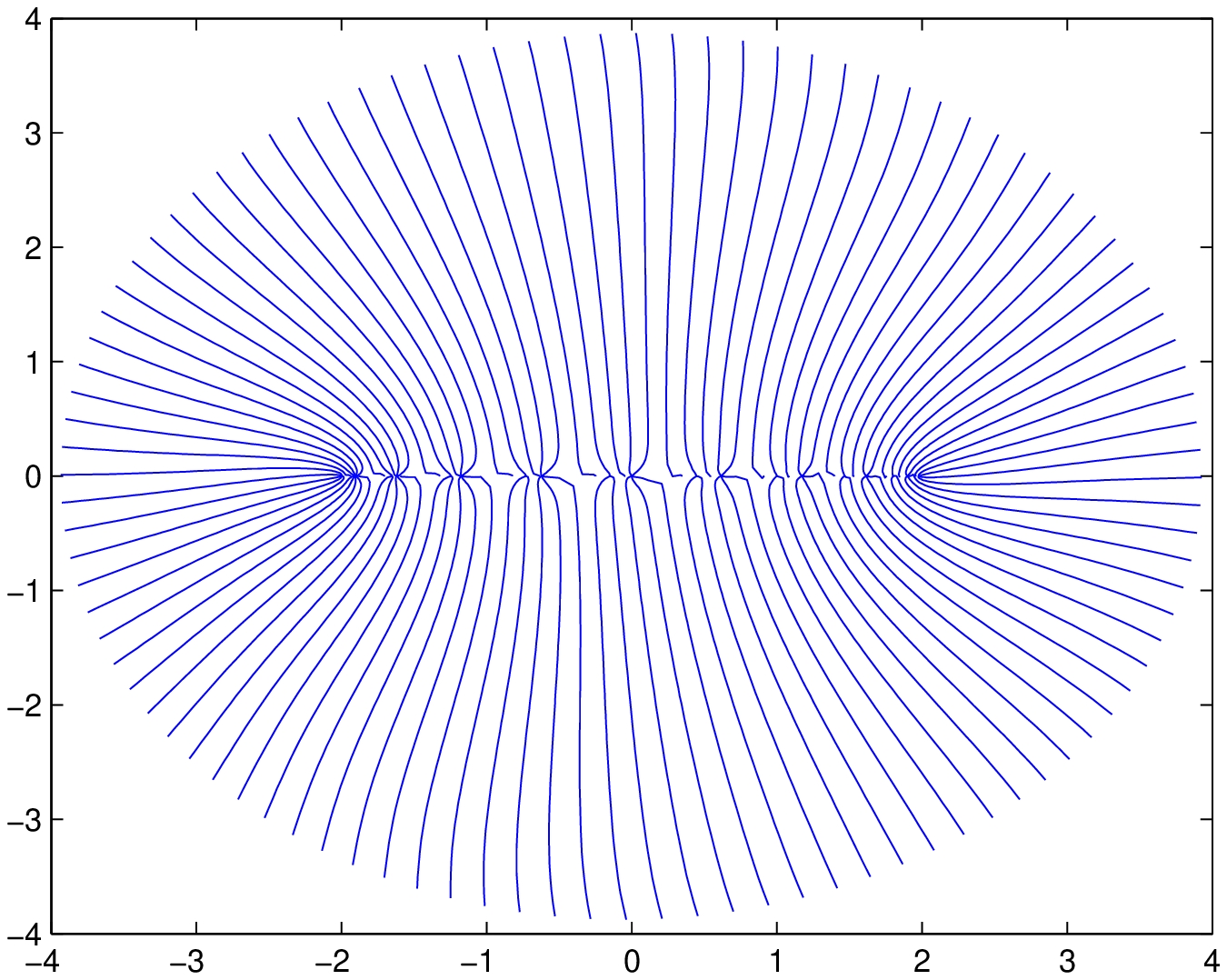}{20}
\fig{1}{History of the convergence for the H problem with $n=5$ and $k=20$}{0}
\end{center}
\medskip

With this choice of starting points we have observed that the number of global  iterations  
is typically   of order of $N$ but   there are   not enough  early deflations, that is, 
iterations that are prematurely stopped due to early convergence.   In order to increase 
the cost  savings due to  premature deflation  in our program  we have employed  a slightly refined strategy.
  Since the method does not approximate directly the eigenvalues $\lambda_i$ but 
their Dickson transform $y_i =\lambda_i +\lambda_i^{-1}$, we have chosen starting points on the Dickson transform of the circles $|z|=\rho$, 
that is points lying on ellipses $\frac{\mathrm{Re}(z)^2}{(\rho+1/\rho)^2} + \frac{\mathrm{Im}(z)^2}{(\rho-1/\rho)^2} = 1$. 
More precisely, this is the algorithm we used to pick the starting points:

\begin{code}
{\bf Input:} Number $N$ of eigenvalues to approximate and parameters $a \in \mathbb N$ and $b\in \mathbb N$\\
{\bf Output:} Starting points $z_k$, $k=1,\dots, N$\\
\> $\theta=2 \pi/N$;\\
\> $\phi$=randn;\\
\> for $j=1,\dots,N$\\
\> \> $jj$=mod($j,a$);\\
\> \> $\rho=1-jj/b$;\\
\> \> $\alpha=\rho+1/\rho$;\\
\> \> $\beta=1/\rho - \rho$;\\
\> \> $z_j=\alpha \cos(j*\theta+\phi)+\beta \sin (j*\theta+\phi)$\\
\> end\\
\end{code}

\vskip-0.6cm

The integer $a$ determines  the number of ellipses  whereas $b$ is used to  
tune the lenghts $\alpha$ and $\beta$, defined as above, of their semiaxes.
 We expect that a good choice for the parameters $a$ and $b$ depends on the ratio $k/n$: when $k \gg n$ we 
expect many eigenvalues to lie on or near to the unit circle, while when $n \gg k$ we expect a situation more similar 
to the eigenvalues of a random matrix, with no particular orientation towards unimodularity. We therefore expect that a small ratio $a/b$ works well in the former case while on the contrary in the latter case $a \simeq b$ should be a better choice. Moreover, we expect that 
as $nk$ grows it is helpful to increase the total number $a$ of ellipses as well. 

We show here some of the results on random  T-palindromic polynomials. Figure 2 refers to an experiment on small-dimensional, high-degree 
polynomials: the value of $n$ has been set to $5$ while $k$ was variable. The average number of $t$ over a set of $1000$ 
random polynomials for each value of $N=nk$ is shown on  the graph. 
The parameters  satisfy  $a\in \{2,3\}$ and $b\in \{8, 64\}$ and  they are determined  by $a=1+2^c$ and $b=8^{c+1}$,  where the integer $c$ 
is defined as $c = \log_{320} N$. 
The graph shows  a linear growth of $t$ with respect to $N=nk$.

\begin{center}
\psone{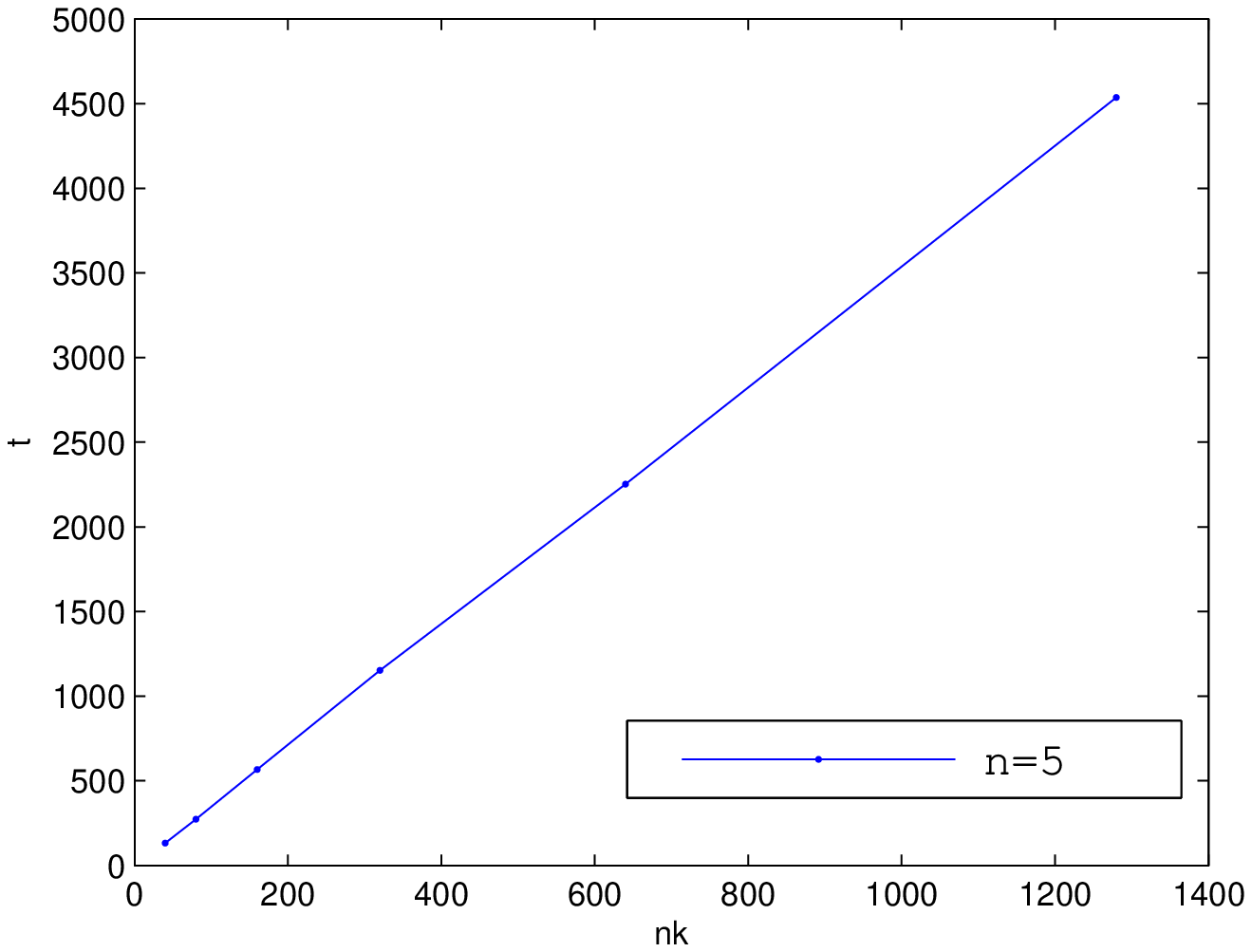}{20}
\end{center}
\fig{2}{}{0}
\medskip

Figure 3 refers to an experiment where on the contrary the case of small $k$  is  explored. 
We  have considered here $k=2$ and let $n$ vary and we show the results for $t$ plotted against $nk$ 
for several choices of $a$ and $b$. The choice labelled as 'step function' is   for $a=\{6, 11\}$ and $b=\{6, 12\}$ generated by 
 $a=1+5 \ 2^c$ and $b=6 \ 2^c$.   Once again the experiments  suggest that when the starting points 
are conveniently chosen $t \leq \alpha N$ for some constant $\alpha$ and any $N$ in the specified range, and,  moreover, 
the  bound  still holds for different 
reasonable  choices of the parameters $a$ and $b$.  The experimentation with random 
polynomials gives $\alpha \simeq 8$ as an estimate for the constant.

\begin{center}
\psone{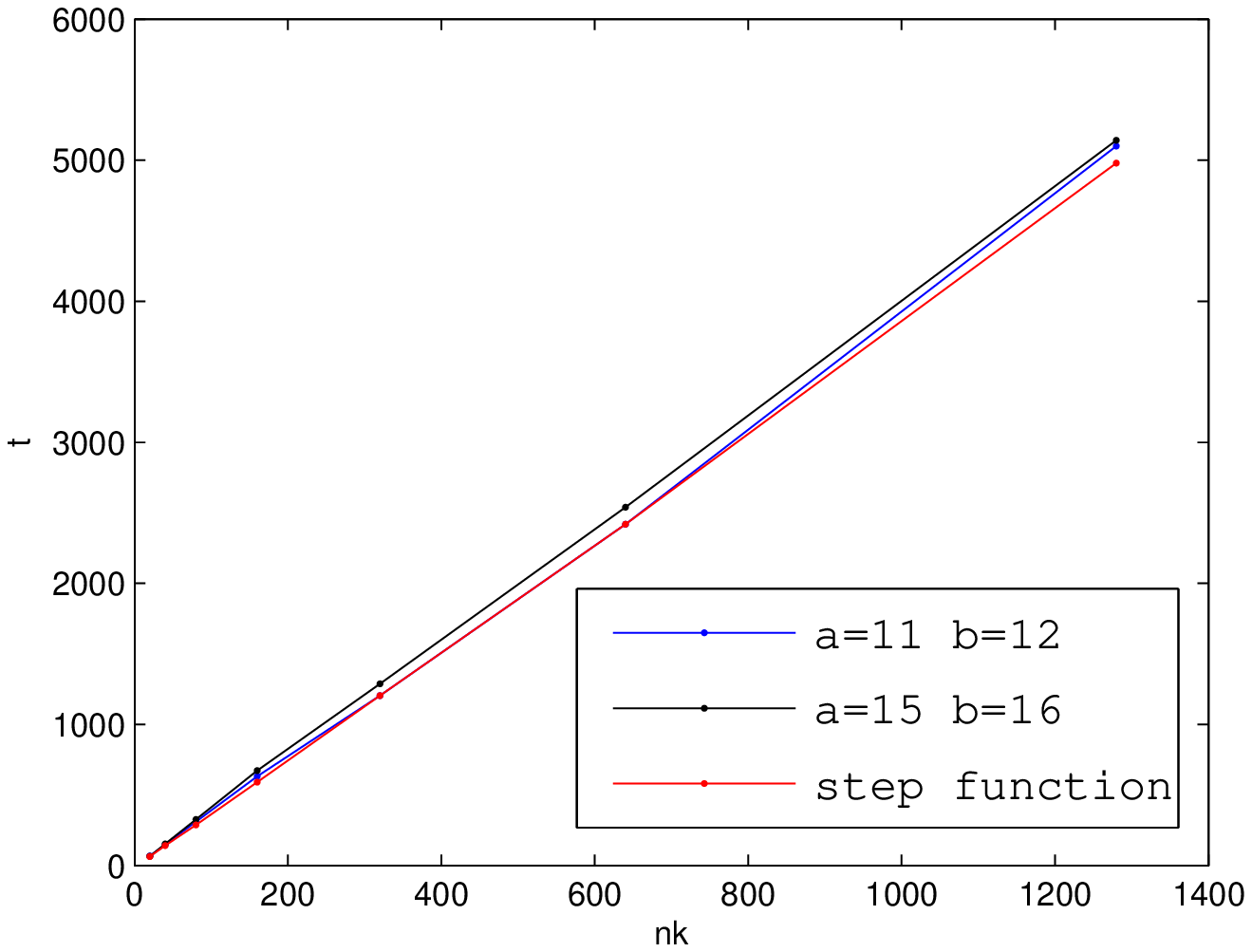}{20}
\end{center}
\fig{3}{}{0}
\medskip

In conclusion,  the algorithm can greatly benefit from a smart  
strategy for the selection of the starting points by increasing the number  of 
early deflations.  The experiments show that as long as 
the starting points are suitably chosen the value of $t$ is proportional to $N=nk$.

\subsection{Accuracy of root-finding}

The other important aspect of our solver based on polynomial root-finding concerns the accuracy of computed approximations.  
 In  our experience the method competes very well in accuracy with the customary 
QZ-algorithm.  The accuracy of the computed non-exceptional roots for the random polynomials was always comparable with the 
accuracy of the approximation obtained with the $QZ$ method. The results of other 
 numerical experiments   confirm 
the  robustness of the novel  method.
Figure 4  illustrates the  computed eigenvalues  for the problem $H_{5, 40}$. 
Figure 5 also reports the  plot of  the absolute error vector ${\emph abs}(\B
\lambda_{EA} -\B {\tilde \lambda})$ and ${\emph abs}(\B
\lambda_{QZ} -\B {\tilde \lambda})$, where $\B {\tilde \lambda}$ is the vector formed by the eigenvalues computed in high precision arithmetic by Mathematica\footnote{Mathematica is a registered trademark of Wolfram Research, Inc.} while $\B \lambda_{EA}$ and 
$\B {\lambda_{QZ}}$ are, respectively, the vectors formed from  the eigenvalues returned  by  our
routine {\bf palindromic$\_$aberth$\_$zeros}  and  suitably sorted by the  internal  function 
${\emph polyeig}$.

\begin{center}
\psone{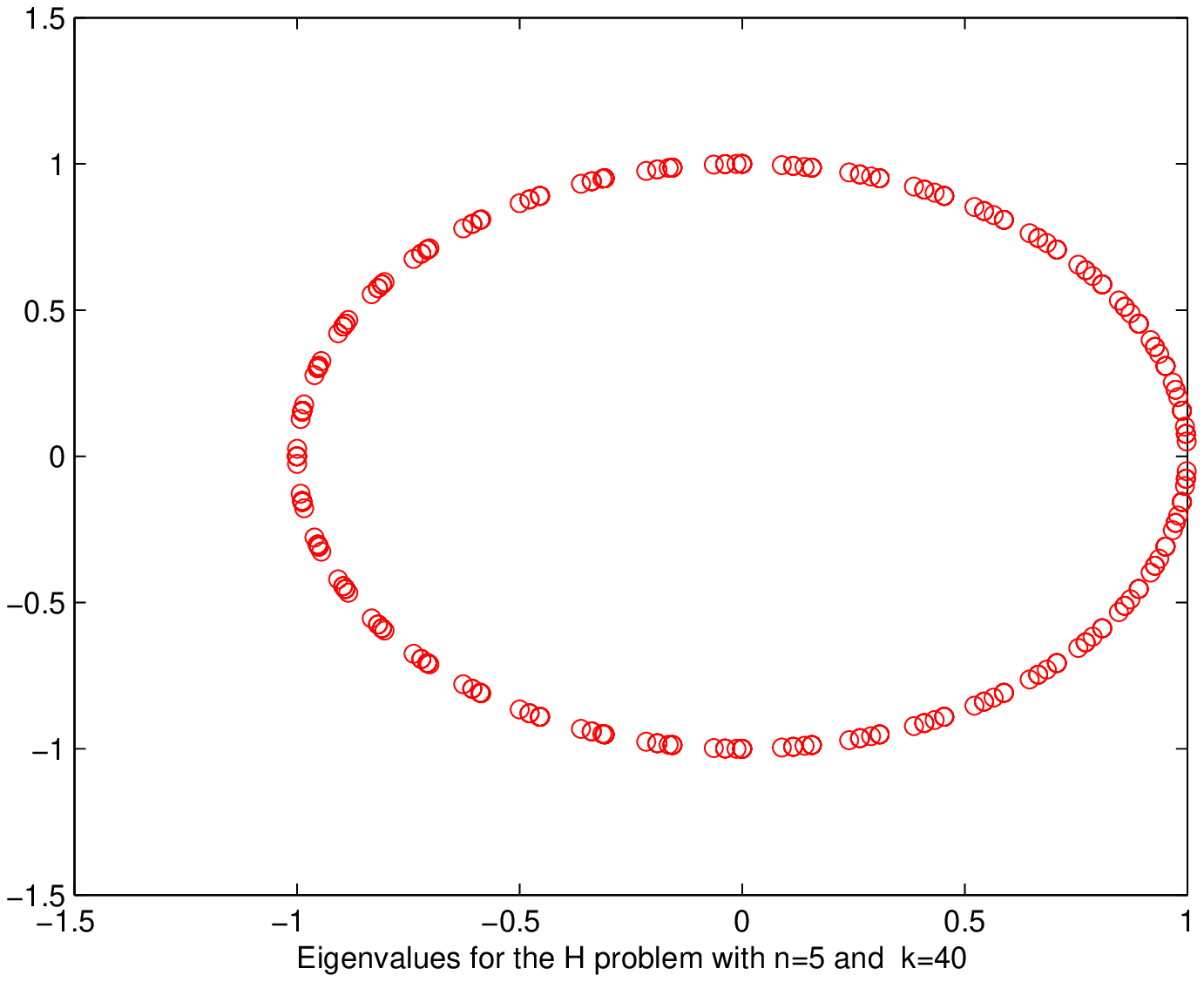}{20}
\end{center}
\fig{4}{}{0}
\medskip

\begin{center}
\psone{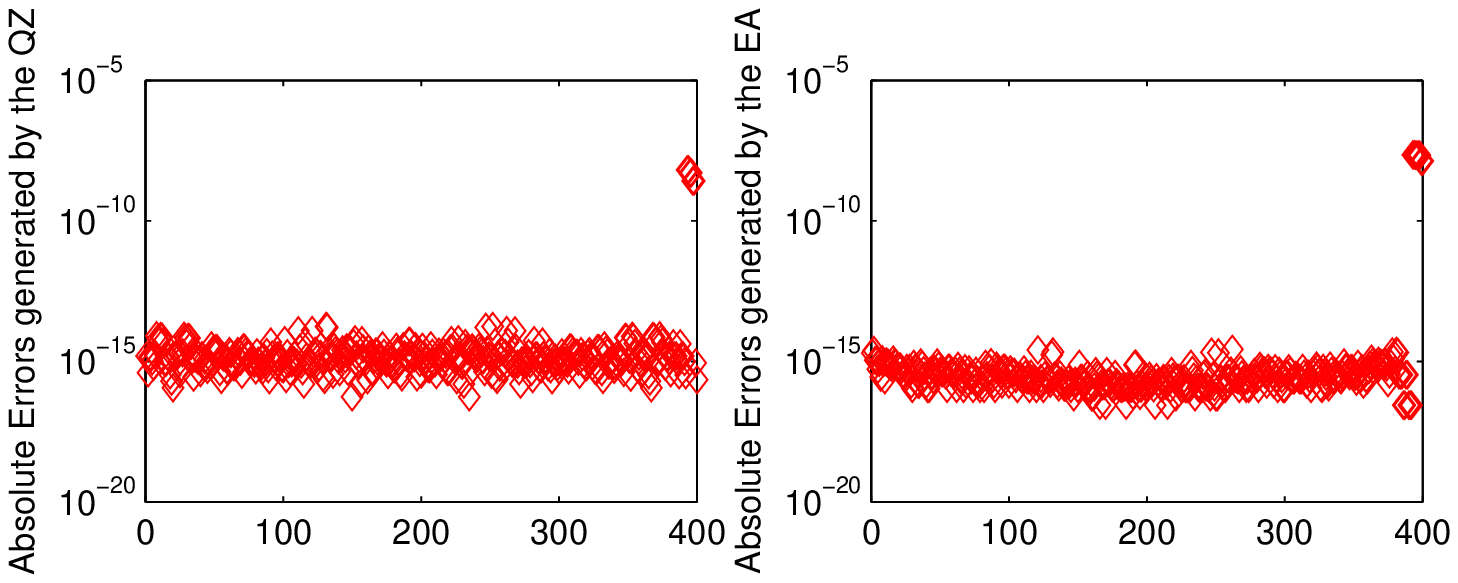}{20}
\end{center}
\fig{5}{}{0}
\medskip

The numerical results  put in evidence the following important  aspects: 
\begin{enumerate}
\item  Poor approximations for the exact eigenvalue
$\lambda=-1$  are in accordance with the  theoretical predictions: 
in fact the reverse transformation from $y=\lambda+\lambda^{-1}$ to $\lambda=\frac{1}{2} (y \pm \sqrt{y^2 - 4})$ 
is known to be ill-conditioned near $y=\pm 2$ (or $\lambda=\pm 1$).   
Since in this example $-1$ is a defective eigenvalue, the approximations returned by ${\emph polyeig}$ have comparable absolute errors of order 
$10^{-8}$ which are  in accordance with the unstructured backward error 
estimates given in \cite{hihi}.

\item  The accuracy  of the remaining  approximations is unaffected  
from the occurrence of  near-to-critical eigenvalues and is in accordance with the 
results  returned by ${\emph polyeig}$. For most non-exceptional eigenvalues, the accuracy of approximations computed by our method is slightly better.
\item This kind of behavior 
is confirmed by many other 
experiments. Our method performs similarly to the QZ for non-exceptional eigenvalues and for defective exceptional eigenvalues, but generally worse than QZ and the structure-preserving methods \cite{urv} for exceptional eigenvalues.

\end{enumerate}

 \section{Conclusions and future work}\label{end}
In this paper we have shown that the  Ehrlich-Aberth  method can be used for 
solving palindromic and 
T-palindromic generalized eigenproblems. The basic idea can be applied to a generic matrix polynomial of any kind; moreover, as we have shown in this paper, it is possible to adapt it in order to exploit certain structures as the palindromic structure that we have considered here. The resulting algorithm  
is numerically robust and  achieves computational  efficiency  
by exploiting the rank-structure of the associated linearization in the Dickson
basis. The algorithm is 
quite interesting for its potential for parallelization on  distributed
architectures and, moreover, can be easily  
incorporated in the MPSolve package to  develop a multiprecision root finder for
matrix polynomial eigenproblems. 

There are, however,  some issues that  still stand in the way of a fully satisfactory  implementation of our method  
and  are currently under  investigation. 
\begin{enumerate}
\item
The development of an automatic procedure for the selection of starting points  is important to 
attain a low operation count due to the prevalence and ease of deflation.  We have shown that a smart 
choice  could  be based on a few parameters to be determined  from some rough information on the 
spectrum localization. 
\item The proposed
algorithm is still
inefficient with respect to the size of the  polynomial coefficients.  The
preliminary reduction of the 
linearized problem into a Hessenberg-triangular form 
is the mean to devise a unified  efficient algorithm for  both small and large
coefficients.   A fast reduction algorithm would be incorporated in our implementation. 
The algorithm  should  be able to exploit the rank structure of the linearization  
(for large degrees), and, at the same time, the inner 
structure of the quasiseparable generators  (for large coefficients).
\item  Regarding the accuracy  of the method there are still some difficulties in the numerical 
treatment of the critical cases.  Our current research is focusing on the issue of a structured refinement of the approximations of such eigenvalues.
\end{enumerate}

\bigskip
\noindent {\bf Acknowledgements.}
The second author wishes  to thank 
Volker Mehrmann for the many discussions and the valuable comments and suggestions 
provided during  the 
Gene Golub SIAM Summer School 2010. We are grateful to the anonymous referees for their hints and suggestions.

\end{document}